\documentclass[reqno,11pt]{amsart}
\usepackage{xcolor}

\usepackage{amssymb}
\usepackage{comment}
\usepackage{tikz}
\usepackage{bbm}
\usepackage[all,cmtip]{xy}
\newcounter{mt}

\newtheorem{Proposition}{Proposition}[section]
\newtheorem{Definition}[Proposition]{Definition}
\newtheorem{Lemma}[Proposition]{Lemma}
\newtheorem{Theorem}[Proposition]{Theorem}
\newtheorem{Corollary}[Proposition]{Corollary}

\newtheorem{Conjecture}[Proposition]{Conjecture}

\DeclareMathOperator{\Gr}{Gr}

\DeclareMathOperator{\id}{Id}

\DeclareMathOperator{\vol}{vol}
\DeclareMathOperator{\Lin}{Lin}

\DeclareMathOperator{\GL}{GL}

\DeclareMathOperator{\OO}{O}
\DeclareMathOperator{\SO}{SO}

\DeclareMathOperator{\cent}{cent}

\newcommand{\R}{\mathbb{R}}

\newcommand{\calE}{\mathcal{E}}

\newcommand{\calK}{\mathcal{K}}

\makeatletter
\def\moverlay{\mathpalette\mov@rlay}
\def\mov@rlay#1#2{\leavevmode\vtop{%
		\baselineskip\z@skip \lineskiplimit-\maxdimen
		\ialign{\hfil$\m@th#1##$\hfil\cr#2\crcr}}}
\newcommand{\charfusion}[3][\mathord]{
	#1{\ifx#1\mathop\vphantom{#2}\fi
		\mathpalette\mov@rlay{#2\cr#3}
	}
	\ifx#1\mathop\expandafter\displaylimits\fi}
\makeatother

\title[Stability in the Banach isometric conjecture]{Stability in the Banach isometric conjecture and nearly monochromatic Finsler surfaces}
\author{Gautam Aishwarya}
	\email{gautamaish@gmail.com} 
\address{Faculty of Mathematics, Technion - Israel Institute of Technology, Haifa 3200003, Israel}
\author{Dmitry Faifman}
	\email{faifmand@tauex.tau.ac.il} 
\address{School of Mathematical Sciences, Tel Aviv University, Tel Aviv 6997801, Israel}

\begin{document}
	
	\begin{abstract}
		The Banach isometric conjecture asserts that a normed space with all of its $k$-dimensional subspaces isometric, where $k\geq 2$, is necessarily Euclidean. The first case of $k=2$ is classical, established by Auerbach, Mazur and Ulam using an elegant topological argument. 
		We refine their method to arrive at a stable version of their result: if all $2$-dimensional subspaces are almost isometric, then the space is almost Euclidean.
		Furthermore, we show that a  $2$-dimensional surface, which is not a torus or a Klein bottle, equipped with a near-monochromatic Finsler metric, is approximately Riemannian. The stability is quantified explicitly using the Banach--Mazur distance. 
	\end{abstract}

\thanks{  
	\\\indent GA was supported by ISF grant 1468/19 and NSF-BSF grant DMS-2247834. DF was supported by ISF grant No. 1750/20. 
}

\maketitle

\section{Introduction and results}
In 1932, Banach made the following conjecture \cite{Banach30}. 
\begin{Conjecture}\label{conj:banach1}
	Let $V$ be a normed space over the reals, and  $2\leq k<\dim V$ an integer. Suppose that all linear $k$-dimensional subspaces of $V$ are isometric to each other. Then $V$ must necessarily be a Hilbert space. 
\end{Conjecture}
A seemingly more general conjecture, known to be equivalent to Conjecture \ref{conj:banach1} \cite{Montejano21arxivsurvey}, is the following.
\begin{Conjecture}\label{conj:banach2}
Let $V$ be a finite-dimensional real vector space, and $2\leq k<\dim V$. If $K\subset V$  is a convex body, and all its $k$-dimensional sections through a fixed interior point of $K$ are affinely equivalent, then $K$ is an ellipsoid.
\end{Conjecture}
Here two subsets of two respective affine spaces are \emph{affinely equivalent} if one is the image of the other under an affine map.

This conjecture has been proven in various cases, starting with the work of Auerbach--Mazur--Ulam in 1935 \cite{AuerbachMazurUlam35}, who proved it for $k=2$. Their proof is based on a topological obstruction - an idea later pushed further by Gromov \cite{Gromov67} to settle the conjecture for all even $k$, as well as for $\dim V\geq k+2$. Much more recently, Bor--Hernandez Lamoneda--Jimenez Desantiago--Montejano Peimbert  \cite{BorLamonedaDesantiagoMontejano21} combined topological and convex-geometric ideas to prove the conjecture for $k=4m+1$, with the possible exception of $k=133$; while for $k=3$, the conjecture has been confirmed by Ivanov, Mamaev, and Nordskova \cite{IvanovMamaevNordskova23}. A local version of Conjecture \ref{conj:banach1} for $k=2,3$ was established in \cite{Ivanov18, ivanov2023local}. Some results were also obtained for the Banach conjecture in complex normed spaces \cite{Gromov67, montejano_complex}.

\subsection{Results}

In this note, we initiate the study of stability in Banach's conjecture. Stability results in convex geometry concern themselves either with geometric inequalities, showing near-extremizers to be close to the extremizers; or with the geometric characterization of a certain type of convex body, showing e.g. that a convex body approximately satisfying a property exhibitied only by ellipsoids, must be close to an ellipsoid.  This work studies the latter type of stability. Previous stability results of this type appeared in several works, such as Burger and Schneider \cite{BurgerSchneider93}, Groemer \cite{Groemer94}, Gruber \cite{Gruber97} (see also references therein).

Loosely speaking, we wish to know if, when the linear $k$-dimensional sections of a convex body $K$ with $0$ in its interior are approximately linearly/affinely equivalent, $K$ must be approximately an ellipsoid. To quantify the deviation from linear/affine equivalence, we use the Banach--Mazur distance.

\begin{Definition}
Let $V,W$ be isomorphic Banach spaces (with possibly non-symmetric norms). The Banach--Mazur distance between them is defined by
\[
d_{BM}(V,W) = \inf \{ \Vert T \Vert \Vert T^{-1} \Vert: T \in \GL (V,W) \}.
\]
\end{Definition}
One can identify (non-symmetric) normed spaces with their unit balls, which are convex bodies with the origin a fixed point in their interior. The Banach-Mazur distance yields a corresponding notion of distance on such pointed bodies.
\begin{Definition} Let $K\subset V$, $L\subset W$ be convex bodies with the origin in their interior. The \emph{linear Banach--Mazur distance} is
	\[
	d^{\Lin}_{BM} (K,L) = \inf \{ \lambda>0 : K \subset TL \subset \lambda K, \textnormal{ for some }T \in {\GL} (W, V). \}.
	\]
\end{Definition}

When affine maps are considered, we need to allow additional freedom in the choice of origin. 
Denote by $\overline{\GL}(W, V)$ the set of affine maps from $W$ to $V$. A convex body in $\R^n$ is any compact convex set with non-empty interior.
\begin{Definition} Let $K\subset V$, $L\subset W$ be convex bodies.
\[
d_{BM} (K,L) = \inf \{ \lambda>0 : K \subset TL \subset \lambda K+z, \textnormal{ for some }T \in \overline{\GL} (W, V), z\in V \}.
\]
\end{Definition}
Note that if the convex bodies $K,L \subset \R^n$ are the unit balls of norms $\|\bullet\|_K$, resp $\|\bullet \|_L$ on $\R^n$, we have $d_{BM}(K, L)\leq d^{\Lin}_{BM}(K, L)=d_{BM}(\|\bullet\|_K, \|\bullet\|_L)$ in general. Equality holds e.g. when both $K, L$ are centrally symmetric.  Note also that the Banach-Mazur distance is only a non-degenerate distance function on the quotient of the space of convex bodies with non-empty interior by $\overline{\GL}(\R^n)$.

The Banach--Mazur distance is usually left in multiplicative form. Thus $K, L$ are affinely equivalent if and only if $d_{BM}(K, L)=1$, and $d_{BM}(K, L)\geq 1$ in general.

We can now formulate a quantitative version of the Banach conjecture. Denote by $\mathbb E^n$ the Euclidean $n$-dimensional space, and by $B^n$ the unit ball therein.
\begin{Conjecture}\label{conj1}
	For any $\epsilon>0$ and $n>k\geq 2$ there is $\delta=\delta(n, k, \epsilon)>0$ such that the following holds. For any $n$-dimensional normed space $V$, possibly non-symmetric, such that $d_{BM}(E,  F)<1+ \delta$ for all $k$-dimensional linear subspaces $E, F\subset V$, it must hold that
	$d_{BM}(V, \mathbb E^n)<1+\epsilon$.
\end{Conjecture}
As before, this is the linear case of the general affine conjecture.
\begin{Conjecture}\label{conj2}
	For any $\epsilon>0$ and $n>k\geq 2$ there is $\delta=\delta(n, k,\epsilon)>0$ such that the following holds. For any convex body $K\subset\R^n$ with $0$ in its  interior, such that $d_{BM}(K\cap E,  K\cap F)<1+ \delta$ for all $k$-dimensional linear subspaces $E, F\subset \R^n$, we must have
	$d_{BM}(K, B^n)<1+\epsilon$.
\end{Conjecture}

We prove Conjecture \ref{conj2} for $k=2$ with an explicit stability estimate. In the following, $c>0$ denotes various universal constants that can be computed explicitly.

\begin{Theorem}\label{thm:main1}
	Let $K\subset\R^n$ be a convex body with $0$ in its interior. Assume that $d_{BM}(K\cap E, K\cap F)<1+ \delta$ for all linear $2$-dimensional subspaces $E, F\subset \R^n$. Then
	$d_{BM}(K, B^n)<1+c n^{2n^2}\delta^{1/6}$.
\end{Theorem}
In the linear setting of Conjecture \ref{conj1} we obtain better stability.

\begin{Theorem}\label{thm:main_norm}
	Let $V=\R^n$ be a normed space, not necessarily symmetric. Assume that $d_{BM}(E, F)<1+ \delta$ for all $2$-dimensional subspaces $E, F\subset V$. Then
	$d_{BM}(V, \mathbb E^n)<1+cn^2\delta^{1/3}$.
\end{Theorem}
In the special case of symmetric norms this can be stated as follows.
\begin{Corollary}\label{cor:symm}
	Let $K\subset\R^n$ be a centrally-symmetric convex body. Assume that $d_{BM}(K\cap E, K\cap F)<1+ \delta$ for all linear $2$-dimensional subspaces $E, F\subset \R^n$. Then $d_{BM}(K, B^n)<1+cn^2\delta^{1/3}$.
\end{Corollary}

We remark that the main interest of these results lies in the regime of fixed $n$ and $\delta\to 0$. Indeed by the Dvoretzky--Milman theorem \cite[Theorem 5.13]{artstein_etal}, if $n>\exp(c\delta^{-2})$ and $K$ is a centrally-symmetric convex body, then $K$ must have a planar section $E_0$ such that $d_{BM}(K\cap E_0, B^2)<1+\delta$; under the assumptions of Corollary \ref{cor:symm}, we have $d_{BM}(K\cap E, B^2)<1+c\delta$ for all planar sections, yielding $d_{BM}(K, B^n)<1+cn^2\delta$ as in the proof of Theorem \ref{thm:main_norm}.

We deduce these results as a corollary of the following general theorem on fields of convex bodies. A \emph{continuous field of convex bodies} over a manifold $\Sigma$ is a family of convex bodies $K_x\subset T_x\Sigma$, $x\in\Sigma$, continuous with respect to the Hausdorff distance within each chart.

\begin{Theorem}\label{thm:main_g}
	Let $\Sigma$ be a closed surface which is not a torus or a Klein bottle. If a continuous field of convex bodies $(K_x)_{x\in \Sigma}$ satisfies $d_{BM}(K_x,K_y)<1+\delta$ for all $x,y\in \Sigma$, then for all $x\in \Sigma$, $d_{BM}(K_x, B^2)<1+c\delta^{1/3}$.
\end{Theorem}
It follows that a near-monochromatic Finsler metric on $\Sigma$ is approximately Riemannian: 
\begin{Corollary}\label{cor:monochromatic}
	Let $\Sigma$ be as in Theorem \ref{thm:main_g}, and assume it is equipped with a $C^0$ Finsler structure such that $d_{BM}(T_x\Sigma, T_y\Sigma)<1+\delta$ for all $x, y\in\Sigma$. Then $d_{BM}(T_x\Sigma, \mathbb E^2)<1+c\delta^{1/3}$ for all $x\in\Sigma$.
\end{Corollary}
The restriction on the topology of $\Sigma$ is necessary as one can easily construct monochromatic Finsler structures on the torus and the Klein bottle which are not Riemannian. It will be seen from the proof that the key required property is the non-existence of a field of tangent lines on the surface.
We conjecture that Theorem \ref{thm:main_g}, with some stability bounds, holds in greater generality, in particular for even dimensional spheres.

Let us remark that the power $\frac13$ of $\delta$ established in Theorem \ref{thm:main_g} is likely an artifact of the proof. It would be interesting if the sharp power is strictly less than $1$.

The main step in deducing Theorem \ref{thm:main1} from Theorem \ref{thm:main_g} consists of proving the following stability result, which may be of independent interest.

\begin{Theorem}\label{mprop:one_center}
	Let $K\subset \R^n$ be a convex body with $0\in\mathrm{int}(K)$, such that for all $2$-dimensional linear planes $E$, $d_{BM}(K\cap E, B^2)<1+\epsilon$.
	Then $d_{BM}(K, B^n)<1+cn^{2n^2}\sqrt \epsilon$.
\end{Theorem}
\subsection{Ideas and plan of the proof}
The proof is inspired by the proof of Auerbach--Mazur--Ulam of Conjecture \ref{conj2} for $k=2$, $\dim V=3$, which we now recall using modern language. 
Setting $K_x=K\cap x^\perp$, we get a continuous field of convex bodies in the tangent planes of $S^2$, $x\mapsto K_x\subset T_xS^2$.
A convex body in $\R^2$ is either an ellipse, or has a finite group of affine symmetries. Fixing $x_0\in S^2$ and assuming that none of $K_x$ are ellipses, the affine maps $\mathrm{Aff}(K_{x_0}, K_x)$ mapping $K_{x_0}$ to $K_x$ define a covering space over $S^2$, which must have a section, $x\mapsto g_x$, as $S^2$ is simply connected. Fixing a line $L\subset T_{x_0}S^2$, and letting $L_x$ be the linear line in $T_xS^2$ parallel to $g_x(L)$, we obtain a field $x\mapsto L_x$ of tangent lines on $S^2$, which is impossible. Thus all sections $K_x$ are ellipses. A separate geometric argument, showing that a convex body all of whose $2$-dimensional sections through an interior point are ellipses must be an ellipsoid, completes the proof.

Since the proof of Auerbach--Mazur--Ulam, as do most other known cases of the Banach conjecture, has a topological obstruction at its heart, it is somewhat surprising that a stability result for $k=2$ can be obtained by building on top of the proof outlined above.  We are forced to work with approximate isometries that do not form a group, and thus we aim to distill certain coarse features that are sufficiently stable as to allow the application of topological methods.

The proof of Theorem \ref{thm:main_g} appears in section \ref{sec:surfaces}. We make use of an associated Riemannian structure on $\Sigma$ given by the Binet--Legendre (inertia) ellipsoids $\calE(K_x)$. We consider sets of approximate isometries $G_\alpha(K_x, K_y)$, consisting of affine maps $g:T_x\Sigma\to T_y\Sigma$ that preserve the Binet--Legendre ellipsoids, and map $K_x$ to $K_y$ up to an error $\alpha$ with respect to the Hausdorff distance. We observe that if the elements of $G_{\alpha}(K_x)$ form a fine net in the orthogonal group, then $K_x$ must be close to an ellipse. Provided this is not the case, the key step is then to find a range of values $\alpha$ for which $G_\alpha(K_{x_0})$, or more precisely a low-resolution version thereof, is independent of $\alpha$. This then allows to define a covering space consisting of approximate isometries in $G_{\alpha'}(K_{x_0}, K_x)$, and proceed as before. 

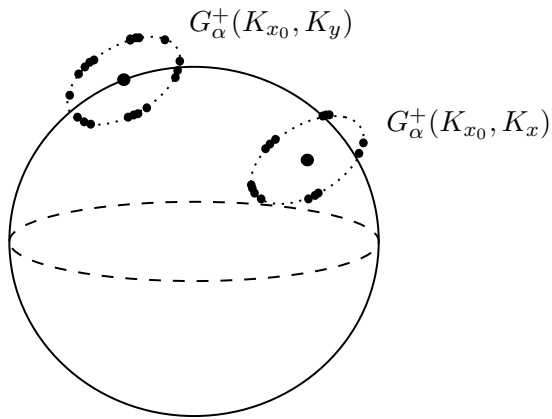
\begin{figure}
\tikzset{every picture/.style={line width=0.75pt}} 

\begin{tikzpicture}[x=0.75pt,y=0.75pt,yscale=-1,xscale=1]

\draw   (357,921.87) .. controls (357,873.19) and (398.7,833.74) .. (450.14,833.74) .. controls (501.59,833.74) and (543.29,873.19) .. (543.29,921.87) .. controls (543.29,970.54) and (501.59,1010) .. (450.14,1010) .. controls (398.7,1010) and (357,970.54) .. (357,921.87) -- cycle ;
\draw  [fill={rgb, 255:red, 0; green, 0; blue, 0 }  ,fill opacity=1 ] (504.43,880.78) .. controls (504.43,879.28) and (505.71,878.07) .. (507.29,878.07) .. controls (508.87,878.07) and (510.14,879.28) .. (510.14,880.78) .. controls (510.14,882.27) and (508.87,883.48) .. (507.29,883.48) .. controls (505.71,883.48) and (504.43,882.27) .. (504.43,880.78) -- cycle ;
\draw  [dash pattern={on 0.84pt off 2.51pt}] (526.16,859.13) .. controls (537.9,861.8) and (538.97,873.65) .. (528.55,885.61) .. controls (518.13,897.56) and (500.16,905.09) .. (488.42,902.42) .. controls (476.67,899.75) and (475.6,887.9) .. (486.02,875.94) .. controls (496.44,863.99) and (514.41,856.46) .. (526.16,859.13) -- cycle ;
\draw  [dash pattern={on 0.84pt off 2.51pt}] (433.58,819.12) .. controls (445.33,821.79) and (446.4,833.64) .. (435.98,845.6) .. controls (425.56,857.55) and (407.59,865.08) .. (395.85,862.41) .. controls (384.1,859.74) and (383.03,847.89) .. (393.45,835.93) .. controls (403.87,823.98) and (421.84,816.45) .. (433.58,819.12) -- cycle ;
\draw  [fill={rgb, 255:red, 0; green, 0; blue, 0 }  ,fill opacity=1 ] (417.57,840.22) .. controls (417.57,838.73) and (416.29,837.52) .. (414.71,837.52) .. controls (413.14,837.52) and (411.86,838.73) .. (411.86,840.22) .. controls (411.86,841.72) and (413.14,842.93) .. (414.71,842.93) .. controls (416.29,842.93) and (417.57,841.72) .. (417.57,840.22) -- cycle ;
\draw  [fill={rgb, 255:red, 0; green, 0; blue, 0 }  ,fill opacity=1 ] (416.34,820.26) .. controls (416.34,819.28) and (417.08,818.49) .. (418,818.49) .. controls (418.92,818.49) and (419.67,819.28) .. (419.67,820.26) .. controls (419.67,821.24) and (418.92,822.03) .. (418,822.03) .. controls (417.08,822.03) and (416.34,821.24) .. (416.34,820.26) -- cycle ;
\draw  [fill={rgb, 255:red, 0; green, 0; blue, 0 }  ,fill opacity=1 ] (420.77,819.08) .. controls (420.77,818.11) and (421.52,817.31) .. (422.44,817.31) .. controls (423.36,817.31) and (424.1,818.11) .. (424.1,819.08) .. controls (424.1,820.06) and (423.36,820.85) .. (422.44,820.85) .. controls (421.52,820.85) and (420.77,820.06) .. (420.77,819.08) -- cycle ;
\draw  [fill={rgb, 255:red, 0; green, 0; blue, 0 }  ,fill opacity=1 ] (419.11,819.08) .. controls (419.11,818.11) and (419.86,817.31) .. (420.77,817.31) .. controls (421.69,817.31) and (422.44,818.11) .. (422.44,819.08) .. controls (422.44,820.06) and (421.69,820.85) .. (420.77,820.85) .. controls (419.86,820.85) and (419.11,820.06) .. (419.11,819.08) -- cycle ;
\draw  [fill={rgb, 255:red, 0; green, 0; blue, 0 }  ,fill opacity=1 ] (416.34,819.67) .. controls (416.34,818.7) and (417.08,817.9) .. (418,817.9) .. controls (418.92,817.9) and (419.67,818.7) .. (419.67,819.67) .. controls (419.67,820.65) and (418.92,821.44) .. (418,821.44) .. controls (417.08,821.44) and (416.34,820.65) .. (416.34,819.67) -- cycle ;
\draw  [fill={rgb, 255:red, 0; green, 0; blue, 0 }  ,fill opacity=1 ] (398.04,830.28) .. controls (398.04,829.3) and (398.78,828.51) .. (399.7,828.51) .. controls (400.62,828.51) and (401.37,829.3) .. (401.37,830.28) .. controls (401.37,831.25) and (400.62,832.05) .. (399.7,832.05) .. controls (398.78,832.05) and (398.04,831.25) .. (398.04,830.28) -- cycle ;
\draw  [fill={rgb, 255:red, 0; green, 0; blue, 0 }  ,fill opacity=1 ] (395.82,831.46) .. controls (395.82,830.48) and (396.56,829.69) .. (397.48,829.69) .. controls (398.4,829.69) and (399.15,830.48) .. (399.15,831.46) .. controls (399.15,832.43) and (398.4,833.22) .. (397.48,833.22) .. controls (396.56,833.22) and (395.82,832.43) .. (395.82,831.46) -- cycle ;
\draw  [fill={rgb, 255:red, 0; green, 0; blue, 0 }  ,fill opacity=1 ] (393.43,833.81) .. controls (393.43,832.84) and (394.17,832.05) .. (395.09,832.05) .. controls (396.01,832.05) and (396.75,832.84) .. (396.75,833.81) .. controls (396.75,834.79) and (396.01,835.58) .. (395.09,835.58) .. controls (394.17,835.58) and (393.43,834.79) .. (393.43,833.81) -- cycle ;
\draw  [fill={rgb, 255:red, 0; green, 0; blue, 0 }  ,fill opacity=1 ] (390.27,837.35) .. controls (390.27,836.37) and (391.02,835.58) .. (391.94,835.58) .. controls (392.86,835.58) and (393.6,836.37) .. (393.6,837.35) .. controls (393.6,838.33) and (392.86,839.12) .. (391.94,839.12) .. controls (391.02,839.12) and (390.27,838.33) .. (390.27,837.35) -- cycle ;
\draw  [fill={rgb, 255:red, 0; green, 0; blue, 0 }  ,fill opacity=1 ] (389.72,859.15) .. controls (389.72,858.18) and (390.46,857.38) .. (391.38,857.38) .. controls (392.3,857.38) and (393.05,858.18) .. (393.05,859.15) .. controls (393.05,860.13) and (392.3,860.92) .. (391.38,860.92) .. controls (390.46,860.92) and (389.72,860.13) .. (389.72,859.15) -- cycle ;
\draw  [fill={rgb, 255:red, 0; green, 0; blue, 0 }  ,fill opacity=1 ] (393.05,861.51) .. controls (393.05,860.53) and (393.79,859.74) .. (394.71,859.74) .. controls (395.63,859.74) and (396.37,860.53) .. (396.37,861.51) .. controls (396.37,862.48) and (395.63,863.28) .. (394.71,863.28) .. controls (393.79,863.28) and (393.05,862.48) .. (393.05,861.51) -- cycle ;
\draw  [fill={rgb, 255:red, 0; green, 0; blue, 0 }  ,fill opacity=1 ] (396.37,862.69) .. controls (396.37,861.71) and (397.12,860.92) .. (398.04,860.92) .. controls (398.96,860.92) and (399.7,861.71) .. (399.7,862.69) .. controls (399.7,863.66) and (398.96,864.45) .. (398.04,864.45) .. controls (397.12,864.45) and (396.37,863.66) .. (396.37,862.69) -- cycle ;
\draw  [fill={rgb, 255:red, 0; green, 0; blue, 0 }  ,fill opacity=1 ] (420.22,857.38) .. controls (420.22,856.41) and (420.97,855.62) .. (421.88,855.62) .. controls (422.8,855.62) and (423.55,856.41) .. (423.55,857.38) .. controls (423.55,858.36) and (422.8,859.15) .. (421.88,859.15) .. controls (420.97,859.15) and (420.22,858.36) .. (420.22,857.38) -- cycle ;
\draw  [fill={rgb, 255:red, 0; green, 0; blue, 0 }  ,fill opacity=1 ] (418,858.05) .. controls (418,857.08) and (418.75,856.29) .. (419.67,856.29) .. controls (420.58,856.29) and (421.33,857.08) .. (421.33,858.05) .. controls (421.33,859.03) and (420.58,859.82) .. (419.67,859.82) .. controls (418.75,859.82) and (418,859.03) .. (418,858.05) -- cycle ;
\draw  [fill={rgb, 255:red, 0; green, 0; blue, 0 }  ,fill opacity=1 ] (415.23,859.15) .. controls (415.23,858.18) and (415.97,857.38) .. (416.89,857.38) .. controls (417.81,857.38) and (418.56,858.18) .. (418.56,859.15) .. controls (418.56,860.13) and (417.81,860.92) .. (416.89,860.92) .. controls (415.97,860.92) and (415.23,860.13) .. (415.23,859.15) -- cycle ;
\draw  [fill={rgb, 255:red, 0; green, 0; blue, 0 }  ,fill opacity=1 ] (424.66,854.44) .. controls (424.66,853.46) and (425.4,852.67) .. (426.32,852.67) .. controls (427.24,852.67) and (427.98,853.46) .. (427.98,854.44) .. controls (427.98,855.41) and (427.24,856.21) .. (426.32,856.21) .. controls (425.4,856.21) and (424.66,855.41) .. (424.66,854.44) -- cycle ;
\draw  [fill={rgb, 255:red, 0; green, 0; blue, 0 }  ,fill opacity=1 ] (440.18,835.58) .. controls (440.18,834.61) and (440.93,833.81) .. (441.85,833.81) .. controls (442.77,833.81) and (443.51,834.61) .. (443.51,835.58) .. controls (443.51,836.56) and (442.77,837.35) .. (441.85,837.35) .. controls (440.93,837.35) and (440.18,836.56) .. (440.18,835.58) -- cycle ;
\draw  [fill={rgb, 255:red, 0; green, 0; blue, 0 }  ,fill opacity=1 ] (441.29,830.87) .. controls (441.29,829.89) and (442.04,829.1) .. (442.96,829.1) .. controls (443.88,829.1) and (444.62,829.89) .. (444.62,830.87) .. controls (444.62,831.84) and (443.88,832.64) .. (442.96,832.64) .. controls (442.04,832.64) and (441.29,831.84) .. (441.29,830.87) -- cycle ;
\draw  [fill={rgb, 255:red, 0; green, 0; blue, 0 }  ,fill opacity=1 ] (439.08,839.12) .. controls (439.08,838.14) and (439.82,837.35) .. (440.74,837.35) .. controls (441.66,837.35) and (442.4,838.14) .. (442.4,839.12) .. controls (442.4,840.09) and (441.66,840.88) .. (440.74,840.88) .. controls (439.82,840.88) and (439.08,840.09) .. (439.08,839.12) -- cycle ;
\draw  [fill={rgb, 255:red, 0; green, 0; blue, 0 }  ,fill opacity=1 ] (516.71,857.97) .. controls (516.71,857) and (517.46,856.21) .. (518.38,856.21) .. controls (519.3,856.21) and (520.04,857) .. (520.04,857.97) .. controls (520.04,858.95) and (519.3,859.74) .. (518.38,859.74) .. controls (517.46,859.74) and (516.71,858.95) .. (516.71,857.97) -- cycle ;
\draw  [fill={rgb, 255:red, 0; green, 0; blue, 0 }  ,fill opacity=1 ] (513.39,858.56) .. controls (513.39,857.59) and (514.13,856.79) .. (515.05,856.79) .. controls (515.97,856.79) and (516.71,857.59) .. (516.71,858.56) .. controls (516.71,859.54) and (515.97,860.33) .. (515.05,860.33) .. controls (514.13,860.33) and (513.39,859.54) .. (513.39,858.56) -- cycle ;
\draw  [fill={rgb, 255:red, 0; green, 0; blue, 0 }  ,fill opacity=1 ] (515.05,858.05) .. controls (515.05,857.08) and (515.8,856.29) .. (516.71,856.29) .. controls (517.63,856.29) and (518.38,857.08) .. (518.38,858.05) .. controls (518.38,859.03) and (517.63,859.82) .. (516.71,859.82) .. controls (515.8,859.82) and (515.05,859.03) .. (515.05,858.05) -- cycle ;
\draw  [fill={rgb, 255:red, 0; green, 0; blue, 0 }  ,fill opacity=1 ] (489.54,870.35) .. controls (489.54,869.37) and (490.29,868.58) .. (491.2,868.58) .. controls (492.12,868.58) and (492.87,869.37) .. (492.87,870.35) .. controls (492.87,871.32) and (492.12,872.12) .. (491.2,872.12) .. controls (490.29,872.12) and (489.54,871.32) .. (489.54,870.35) -- cycle ;
\draw  [fill={rgb, 255:red, 0; green, 0; blue, 0 }  ,fill opacity=1 ] (486.77,872.7) .. controls (486.77,871.73) and (487.51,870.94) .. (488.43,870.94) .. controls (489.35,870.94) and (490.1,871.73) .. (490.1,872.7) .. controls (490.1,873.68) and (489.35,874.47) .. (488.43,874.47) .. controls (487.51,874.47) and (486.77,873.68) .. (486.77,872.7) -- cycle ;
\draw  [fill={rgb, 255:red, 0; green, 0; blue, 0 }  ,fill opacity=1 ] (484.55,875.06) .. controls (484.55,874.08) and (485.29,873.29) .. (486.21,873.29) .. controls (487.13,873.29) and (487.88,874.08) .. (487.88,875.06) .. controls (487.88,876.04) and (487.13,876.83) .. (486.21,876.83) .. controls (485.29,876.83) and (484.55,876.04) .. (484.55,875.06) -- cycle ;
\draw  [fill={rgb, 255:red, 0; green, 0; blue, 0 }  ,fill opacity=1 ] (477.34,893.33) .. controls (477.34,892.35) and (478.09,891.56) .. (479,891.56) .. controls (479.92,891.56) and (480.67,892.35) .. (480.67,893.33) .. controls (480.67,894.3) and (479.92,895.1) .. (479,895.1) .. controls (478.09,895.1) and (477.34,894.3) .. (477.34,893.33) -- cycle ;
\draw  [fill={rgb, 255:red, 0; green, 0; blue, 0 }  ,fill opacity=1 ] (477.89,895.1) .. controls (477.89,894.12) and (478.64,893.33) .. (479.56,893.33) .. controls (480.48,893.33) and (481.22,894.12) .. (481.22,895.1) .. controls (481.22,896.07) and (480.48,896.86) .. (479.56,896.86) .. controls (478.64,896.86) and (477.89,896.07) .. (477.89,895.1) -- cycle ;
\draw  [fill={rgb, 255:red, 0; green, 0; blue, 0 }  ,fill opacity=1 ] (479,897.45) .. controls (479,896.48) and (479.75,895.69) .. (480.67,895.69) .. controls (481.59,895.69) and (482.33,896.48) .. (482.33,897.45) .. controls (482.33,898.43) and (481.59,899.22) .. (480.67,899.22) .. controls (479.75,899.22) and (479,898.43) .. (479,897.45) -- cycle ;
\draw  [fill={rgb, 255:red, 0; green, 0; blue, 0 }  ,fill opacity=1 ] (482.33,900.4) .. controls (482.33,899.42) and (483.08,898.63) .. (484,898.63) .. controls (484.91,898.63) and (485.66,899.42) .. (485.66,900.4) .. controls (485.66,901.38) and (484.91,902.17) .. (484,902.17) .. controls (483.08,902.17) and (482.33,901.38) .. (482.33,900.4) -- cycle ;
\draw  [fill={rgb, 255:red, 0; green, 0; blue, 0 }  ,fill opacity=1 ] (506.18,900.4) .. controls (506.18,899.42) and (506.92,898.63) .. (507.84,898.63) .. controls (508.76,898.63) and (509.5,899.42) .. (509.5,900.4) .. controls (509.5,901.38) and (508.76,902.17) .. (507.84,902.17) .. controls (506.92,902.17) and (506.18,901.38) .. (506.18,900.4) -- cycle ;
\draw  [fill={rgb, 255:red, 0; green, 0; blue, 0 }  ,fill opacity=1 ] (509.5,898.63) .. controls (509.5,897.66) and (510.25,896.86) .. (511.17,896.86) .. controls (512.09,896.86) and (512.83,897.66) .. (512.83,898.63) .. controls (512.83,899.61) and (512.09,900.4) .. (511.17,900.4) .. controls (510.25,900.4) and (509.5,899.61) .. (509.5,898.63) -- cycle ;
\draw  [fill={rgb, 255:red, 0; green, 0; blue, 0 }  ,fill opacity=1 ] (511.17,897.45) .. controls (511.17,896.48) and (511.91,895.69) .. (512.83,895.69) .. controls (513.75,895.69) and (514.5,896.48) .. (514.5,897.45) .. controls (514.5,898.43) and (513.75,899.22) .. (512.83,899.22) .. controls (511.91,899.22) and (511.17,898.43) .. (511.17,897.45) -- cycle ;
\draw  [fill={rgb, 255:red, 0; green, 0; blue, 0 }  ,fill opacity=1 ] (531.69,877.42) .. controls (531.69,876.44) and (532.43,875.65) .. (533.35,875.65) .. controls (534.27,875.65) and (535.01,876.44) .. (535.01,877.42) .. controls (535.01,878.39) and (534.27,879.19) .. (533.35,879.19) .. controls (532.43,879.19) and (531.69,878.39) .. (531.69,877.42) -- cycle ;
\draw  [fill={rgb, 255:red, 0; green, 0; blue, 0 }  ,fill opacity=1 ] (533.91,872.12) .. controls (533.91,871.14) and (534.65,870.35) .. (535.57,870.35) .. controls (536.49,870.35) and (537.23,871.14) .. (537.23,872.12) .. controls (537.23,873.09) and (536.49,873.88) .. (535.57,873.88) .. controls (534.65,873.88) and (533.91,873.09) .. (533.91,872.12) -- cycle ;
\draw  [dash pattern={on 4.5pt off 4.5pt}] (357,921.87) .. controls (357,910.82) and (398.7,901.87) .. (450.14,901.87) .. controls (501.59,901.87) and (543.29,910.82) .. (543.29,921.87) .. controls (543.29,932.91) and (501.59,941.87) .. (450.14,941.87) .. controls (398.7,941.87) and (357,932.91) .. (357,921.87) -- cycle ;
\draw  [color={rgb, 255:red, 0; green, 0; blue, 0 }  ,draw opacity=1 ][fill={rgb, 255:red, 0; green, 0; blue, 0 }  ,fill opacity=1 ] (433.78,820.08) .. controls (433.78,819.11) and (434.53,818.31) .. (435.45,818.31) .. controls (436.37,818.31) and (437.11,819.11) .. (437.11,820.08) .. controls (437.11,821.06) and (436.37,821.85) .. (435.45,821.85) .. controls (434.53,821.85) and (433.78,821.06) .. (433.78,820.08) -- cycle ;
\draw  [color={rgb, 255:red, 0; green, 0; blue, 0 }  ,draw opacity=1 ][fill={rgb, 255:red, 0; green, 0; blue, 0 }  ,fill opacity=1 ] (385.78,848.08) .. controls (385.78,847.11) and (386.53,846.31) .. (387.45,846.31) .. controls (388.37,846.31) and (389.11,847.11) .. (389.11,848.08) .. controls (389.11,849.06) and (388.37,849.85) .. (387.45,849.85) .. controls (386.53,849.85) and (385.78,849.06) .. (385.78,848.08) -- cycle ;

\draw (545.71,853.08) node [anchor=north west][inner sep=0.75pt]    {$G^{+}_{\alpha}( K_{x_{0}} , K_{x})$};
\draw (446.19,802.43) node [anchor=north west][inner sep=0.75pt]    {$G^{+}_{\alpha}( K_{x_{0}} , K_{y})$};

\end{tikzpicture}
\caption{(Oriented) approximate isometries}
\end{figure}

The corollaries of Theorem \ref{thm:main_g} to stability in the Banach conjecture are proved in section \ref{sec:banach}. To pass from Theorem \ref{thm:main_g} to Theorem \ref{thm:main1}, we prove a stable version of independent interest of the characterization of ellipsoids through planar sections mentioned above. For Theorem \ref{thm:main_norm}, we get better stability by utilizing instead a stable version of the von Neumann--Jordan theorem due to Passer.

\subsection{Acknowledgements}
The authors are indebted to Juan Carlos Alvarez Paiva for his inspiring talk at the Integral and Metric Geometry 2022 BIRS-CMO workshop. The second-named author is grateful to Vitali D. Milman for introducing him to the Banach conjecture, and to K\'aroly J. B\"or\"oczky for a helpful discussion.

\section{Preliminaries}
By a convex body in $\R^n$ we understand a compact convex set with non-empty interior. For an ellipsoid $\calE$, we denote by $d_\calE$ the Hausdorff distance with respect to the Euclidean metric for which $\calE$ is a unit ball. Denote by $\cent(K)$ the centroid (center of mass) of $K$. We say $K$ is \emph{centered} if $\cent(K)=0$. A basic property of centered convex bodies in $\R^n$, which we will use repeatedly, is that $-K\subset nK$ \cite[Lemma 2.3.3]{Schneider14}.
For a centered convex body $K$, the \emph{Binet--Legendre ellipsoid} $\calE(K)\subset \R^n$ of $K\subset\R^n$, also called the elipsoid of inertia, is the dual ellipsoid of the unit ball in $(\R^n)^*$ of the inner product
\[\langle \phi, \psi\rangle_K=\frac{n+2}{\vol(K)}\int_K \phi(x)\psi(x)d\vol(x),\]
where $\vol$ is an arbitrary Lebesgue measure on $V$. For general $K$, we set $\calE(K):=\calE(K-\cent(K))+\cent(K)$. 

Evidently, if $g \in \overline{\GL}(\R^n)$ then $\mathcal{E}(gK) = g \mathcal{E}(K)$.
If $K$ is an ellipsoid, $\calE(K)=K$. 
Note that unlike the definition in \cite{MatveevTroyanov12}, which uses a fixed origin, we use the centroid of $K$ as the center of the Binet--Legendre ellipsoid.

\begin{Lemma}\label{lem:BL_continuity}
	Let $K, L\subset\R^n$ be centered convex bodies. If $\frac{1}{\lambda_1}K\subset L\subset  \lambda_2 K$, then 
	\[\calE(K)\subset \lambda_1^{\frac{n+2}{2}}\lambda_2^{\frac n 2}\calE(L),\quad \calE(L)\subset \lambda_1^{\frac n 2}\lambda_2^{\frac {n+2}{2}}\calE(K).\]
	\\Furthermore, there are positive constants $a_n, b_n$ such that $ a_n\calE(K)\subset  K\subset b_n\calE(K)$. 
	One can take explicitly $a_n=2^{-\frac n 2-1}n^{-\frac {3n}4-\frac 32}$, $b_n=2^{\frac n 2+1}n^{\frac {3n}4+ \frac 32}$.
\end{Lemma}
\proof
We follow \cite{MatveevTroyanov15}.
We have
\[\int_{L}\phi(x)^2dx\leq \int_{\lambda_2 K}\phi(x)^2dx=\lambda_2^{n+2}\int_K\phi(x)^2dx. \]
As $\vol(L)\geq \lambda_1^{-n}\vol (K)$, we conclude that
\[\langle \phi, \phi\rangle_L\leq \lambda_1^n\lambda_2^{n+2}\langle \phi, \phi\rangle_K\Rightarrow \calE(L)\subset \lambda_1^{\frac n2}\lambda_2^{\frac {n+2}2}\calE(K).\]
The second inclusion follows by symmetry.

For the second part, define $K_s=\frac12(K+(-K))$. As $-K\subset nK$, we have 
\[   \frac 1n K_s\subset\frac{2}{n+1}K_s \subset K\subset 2K_s,\]
and by the first part,
 \[\calE(K_s)\subset 2^{\frac n2}n^{\frac n 2+1}\calE(K),\quad \calE(K)\subset 2^{\frac n2 +1}n^{\frac n2}\calE(K_s).\]
We make use of John's maximal volume ellipsoid $J$ of $K_s$, which satisfies $J\subset K_s\subset \sqrt nJ$ by John's theorem. Observe that $\calE(J)=J$.
By the first part, it holds that
\[ J\subset n^{\frac n 4}\calE(K_s), \quad \calE(K_s)\subset n^{\frac{n+2}{4}}J.\]
As
\[\frac 1n J\subset \frac 1n K_s\subset K\subset 2K_s\subset 2\sqrt n J,\]
we conclude that
\[K\subset 2^{\frac n 2+1}n^{\frac {3n}4+ \frac 32}\calE(K),\quad \calE(K)\subset 2^{\frac n 2+1}n^{\frac {3n}4+ \frac 32}K.\]
\endproof
\begin{Definition} \label{def:BLdistance}
	The Binet--Legendre distance between $K, L\subset \R^n$ is 
	\[d_{BL}(K, L)=\min \{d_{\calE(g_1K)}(g_1K, g_2L): g_1,g_2\in\overline{\GL}(\R^n), \calE(g_1K)=\calE(g_2L)\}.\]
\end{Definition} 

\begin{Lemma}\label{lem:equivalent_norms}
	$d_{BM}$ and $d_{BL}$ are equivalent distances on $\calK(\R^n)$. Explicitly, there is $\epsilon_n>0$ such that for $ d_{BM}(K, L)\leq1+\epsilon_n$ we have
	\[   a_n(d_{BM}(K, L)-1)\leq d_{BL}(K, L)\leq \widehat b_n(d_{BM}(K, L)-1),\]
	where $\widehat b_n=(3n+1)(1+5n\frac{b_n}{a_n})b_n$.
	
	Furthermore, for $\epsilon<\epsilon_n$, if $d_{BM}(K, L)= 1+\epsilon$ and $\mathcal E(K)=\mathcal E(L)=B^n$, then there is $T\in\OO(\R^n)$ such that 
	\[(1+d_n\epsilon)^{-1} K\subset TL\subset (1+d_n\epsilon) K,\]
	where $d_n=4n\frac{b_n}{a_n}(1+5n\frac{b_n}{a_n})$. 
	
\end{Lemma} 
Later we will use the larger value $d_n=c4^nn^{3n+8}$ for simplicity.
\proof
Write $a=a_n, b=b_n$. Assume $d_{BL}(K, L)=\epsilon$. We may assume that a Euclidean structure is fixed such that $\calE(K)=\calE(L)=B^n$, and $d_{B^n}(K, L)=\epsilon$. That is, $K\subset L+\epsilon B^n$, $L\subset K+\epsilon B^n$.
Then by Lemma \ref{lem:BL_continuity}, $K\subset L+ \frac \epsilon a L=(1+\frac \epsilon a) L$, and similarly $L\subset (1+\frac \epsilon a) K$.
Thus $d_{BM}(K, L)\leq 1+\frac\epsilon a=1+\frac{d_{BL}(K, L)}{a}$.

In the  other direction, assume $d_{BM}(K, L)=1+\epsilon$. We may assume $\calE(K)=B^n$, and $K\subset L\subset (1+\epsilon)K+z$ for some $z\in\R^n$.
It holds that 
\[ q:=\cent(L)=\frac{1}{\vol(L)}\int_{L}xd\vol(x)=\frac{1}{\vol(L)}\int_{L\setminus K}xd\vol(x).\]
 Since $K\subset (1+\epsilon)K+z$, we deduce that 
\[-z\in \epsilon K\Rightarrow z\in n\epsilon K.\]
In particular, $|z|\leq n b\epsilon$, and $\max_L |x|\leq |z|+(1+\epsilon)\max_K |x|\leq b(1+(n+1) \epsilon)$.

Consequently by Lemma \ref{lem:BL_continuity}, and since for $\epsilon<\frac 1{n}$ 
\[(1+\epsilon)^n-1=\sum_{j=1}^n\frac { n(n-1)\cdots(n-j+1)}{j!}\epsilon^j\leq \sum_{j=1}^\infty \frac{(n\epsilon)^j}{j!} =e^{n\epsilon}-1<2 n\epsilon,\]
it follows that for $\epsilon<\frac1{n+1}$,
\[ |q|\leq \frac{ \vol(L\setminus K) }{\vol(L)}\max_{L\setminus K}|x|\leq \frac{((1+\epsilon)^n-1)\vol(K)}{\vol (K)}\max_{L}|x|\leq 4n b \epsilon\]
and so $\pm q\in \frac{4nb}{a}\epsilon aB^n\subset \frac{4nb}{a}\epsilon K$.
Thus 
\[(1-4\frac{nb}{a} \epsilon)K\subset L-q\subset (1+\epsilon)K+\frac{nb\epsilon}{a}aB^n+\frac{4nb\epsilon}{a}aB^n\subset \left(1+ (1+\frac{5nb}{a})\epsilon \right)K.\] 

Denote $b_n'=1+5\frac{nb}{a}$. Then (for $\epsilon$ sufficiently small)
\[ \frac{1}{1+b_n'\epsilon} K\subset L-q\subset (1+b_n'\epsilon)K,\]
and it follows again by Lemma \ref{lem:BL_continuity} that $(1+b_n'\epsilon)^{-n-1}B^n\subset \calE(L-q)\subset (1+b_n'\epsilon)^{n+1}B^n$, and so we can choose a positive definite map $T$ with $T(\calE(L-q))=B^n$ with the eigenvalues of $T$ satisfying $(1+b_n'\epsilon)^{-n-1}\leq\lambda_1, \dots,\lambda_n\leq(1+b_n'\epsilon)^{n+1}$. Writing $T=\id+S$, it follows that the Euclidean operator norm $\|S\|\leq 2nb_n'\epsilon$ (for $\epsilon$ small enough). 
We similarly have $T^{-1}=\id+S'$ with $\|S'\|\leq 2nb_n'\epsilon$.

We deduce that $L'=T(L-q)$ has $\calE(L')=B^n$, and
\begin{align*}L'&\subset T(1+b_n'\epsilon)K=K+b_n'\epsilon K +S(1+b_n'\epsilon) K
	\\&\subset K+ b_n'\epsilon b B^n+b(1+b_n'\epsilon)S(B^n)
	\\&\subset K+ (b_n'\epsilon b +2nb_n'\epsilon b(1+b_n'\epsilon))B^n.\end{align*}
	
	Similarly 
	
	\begin{align*}K&\subset (1+b_n'\epsilon)T^{-1}L'
		\\&\subset L'+ (b_n'\epsilon b +2nb_n'\epsilon b(1+b_n'\epsilon))B^n.\end{align*}
that is, $d_{BL}(K, L)\leq d_{B^n}(K, L')\leq (3n+1)b_n'b\epsilon$. 

In particular, since
\[L'\subset (1+\frac{b}{a}(b_n'\epsilon +2nb_n'\epsilon (1+b_n'\epsilon)))K,\quad K\subset  (1+\frac{b}{a}(b_n'\epsilon +2nb_n'\epsilon (1+b_n'\epsilon)))L',\]
setting $d_n:=4nb_n'\frac{b}{a}$ establishes the last claim of the lemma. 
\endproof

\begin{Corollary}\label{cor:linear_BM}
	If $d_{BM}(K, B^n)<1+\epsilon$ and $\cent(K)=0$ then $(1+d_n\epsilon)^{-1} \mathcal{E}(K) \subset K \subset (1+d_n\epsilon) \mathcal{E}(K)$, and consequently $d_{BM}^{\Lin}(K, B^n)<1+d_n\epsilon$. 
\end{Corollary}

Given an ellipse $\calE\subset \R^2$, we equip $\partial \calE$ with the angular distance function $d_{\calE}$ induced by the Euclidean structure for which it is a unit ball.
Given two ellipses $\calE_i\subset E_i$, $i=1,2$, where $E_i$ are $2$-dimensional linear spaces, denote by $\OO(\calE_1, \calE_2)$ the set of affine maps  $E_1\to E_2$ mapping $\calE_1$ to $\calE_2$. If $E_1=E_2=E$, $\calE_1=\calE_2=\calE$, we write $\OO(\calE)=\OO(\calE, \calE)$, and define $\SO(\calE)\subset\OO(\calE)$ to be the orientation-preserving subgroup.
Note that $\SO(\calE)$ is a circle with distance function $d(g, h)=d_{\calE}(gx,hx)$, which is independent of the choice of $x\in\partial \calE$.
For general $\calE_1, \calE_2$, $\OO(\calE_1, \calE_2)$ is the disjoint union of two circles, and we define a metric on it by
\[  d(g_1,g_2)=\left\{\begin{array}{cc} d_{\calE_1}(g_1^{-1}g_2, \id),&\quad g_1^{-1}g_2\in \SO(\calE_1),\\
											\infty,&\quad  g_1^{-1}g_2\notin \SO(\calE_1).\end{array}\right. \]
Our main technical tool will be the subset of $\lambda$-approximate isometries of convex bodies, defined as follows.
\begin{Definition}
	Let $\lambda>0$, and let $E, F$ be linear spaces of equal dimension. For two convex bodies $K\subset E$, $L\subset F$ we set $G_\lambda(K, L)=\{g\in \OO(\calE(K), \calE(L)): d_{\calE(L)}(L, gK)<\lambda\}$.
	
	We write $G_\lambda(K):=G_\lambda(K, K)\subset\OO(\calE(K))$. We can then also define the subset $G^+_\lambda(K)=G_\lambda(K) \cap \SO(\calE(K))$.
\end{Definition}

The following is immediate from definitions.
\begin{Lemma}
	 If $g\in G_\lambda(K_1, K_2)$ then $g^{-1}\in G_\lambda(K_2, K_1)$. If $g_j\in G_{\lambda_j}(K_j, K_{j+1})$ for $j=1,2$ then $g_2g_1\in G_{\lambda_1+\lambda_2}(K_1, K_3)$.
\end{Lemma}

We say that a set $A\subset S^1$ is an \emph{$\epsilon$-net} if any closed $\epsilon$-interval intersects $A$. 

The following technical lemma is a stable version of the basic fact that a convex body in $\R^2$ is either an ellipse, or has a finite group of symmetries.
\begin{Proposition}\label{prop:dense_euclidean} For all $\epsilon>0$, if $\alpha, \beta>0$ satisfy $\frac 1 {a_2}\alpha+\frac{b_2}{a_2}\beta\leq \epsilon$  and  $G^+_\alpha(K)$ is a $\beta$-net in $\SO(\calE(K))$ then $d_{BM}(K, B^2)<1+\epsilon$.
\end{Proposition}
We will later use \begin{equation}\label{eq:alpha_beta}\alpha(\epsilon)=\frac {a_2}2 \epsilon,\quad \beta(\epsilon)=\frac{a_2}{2b_2}\epsilon,\end{equation} so that $\alpha=\alpha(\epsilon)$, $\beta=\beta(\epsilon)$ satisfy the assumption of the proposition.
\proof
Fix a Euclidean structure $|\bullet|$ on $\R^2$ whose unit ball is $\calE(K)$. Assume $p\in K$ satisfies $|p|=\max\{|x|: x\in K\}$. By Lemma \ref{lem:BL_continuity}, $|p|\leq b_2$. Denoting $D=|p|\calE(K)$, we have $K\subset D$. On the other hand, since $G^+_\alpha(K)$ is a $\beta$-net, it holds that \[D\subset K+2|p|\sin(\frac\beta 2) \calE(K) +\alpha \calE(K) \subset \left( 1+\frac{2b_2}{a_2}\sin\frac \beta2+\frac{1}{a_2}\alpha \right) K,\]
so that $d_{BM}(K, B^2)=d_{BM}(K, D)< 1+\frac{1}{a_2}\alpha+\frac{b_2}{a_2}\beta$.
\endproof

To pass from Theorem \ref{thm:main_g} to Theorem \ref{thm:main_norm}, we make use of the von Neumann--Jordan constant, which we now recall.

\begin{Definition}
	The \emph{von Neumann--Jordan constant} $C_{NJ}(V)$ of a normed space $V$ is the least $M\geq 1$ such that for all non-zero $x, y\in V$,
	\[ \frac{1}{M}\leq  \frac{\|x+y\|^2+\|x-y\|^2}{2(\|x\|^2+\|y\|^2)} \leq M.\]
\end{Definition}
The celebrated von Neumann--Jordan theorem \cite{neumann_jordan} asserts that $C_{NJ}(V)=1$ implies $V$ is Euclidean. The opposite direction is trivial, asserting simply the validity of the parallelogram identity in Euclidean space.

We will need a stable version of these facts, comparing the von Neumann--Jordan constant with the Banach--Mazur distance to Euclidean space. The fact that near-Euclidean normed spaces satisfy an approximate parallelogram identity is easy, and appeared e.g. in \cite[Theorem 5]{KatoMalingradaTakahashi01}.
\begin{Lemma}\label{lem:katoetal} For any normed space $V$ it holds that
	$C_{NJ}(V)\leq d_{BM}(V,\mathbb E^n)^2$.
\end{Lemma}
A stable version of the von Neumann--Jordan theorem was established by Passer. 
\begin{Theorem}[Passer \cite{Passer15}, Theorem 1.2]\label{thm:passer}
For an $n$-dimensional, symmetric normed space $V$ it holds that $d_{BM}(V,\mathbb E^n)\leq 1+cn^2 (C_{NJ}(V)-1)$.
	\end{Theorem}
	Here as before $c>0$ is an explicit universal constant.

\section{Fields of convex bodies on a surface}\label{sec:surfaces}
We now prove Theorem \ref{thm:main_g}. It will be more convenient to work with the following equivalent formulation.
\begin{Theorem}\label{thm:main2'}
	Let $\Sigma$ be a smooth closed surface which is not a torus or a Klein bottle. For any $\delta>0$ and any continuous field of convex bodies $(K_x)_{x\in \Sigma}$, $K_x\subset T_x\Sigma$, such that $d_{BM}(K_x,K_y)<1+\delta$ for all $x,y\in \Sigma$, there must exist $x\in \Sigma$ such that $d_{BM}(K_x, B^2)<1+c\delta^{1/3}$.
\end{Theorem}
	Here and in the proof, $c, C$ are explicit universal constants which assume different values in distinct appearances, whose values we omit for simplicity.
	\proof
	We denote by $d$ the distance function on various circles appearing in the proof. 

	Let $(K_x)$ be a field of convex bodies over $\Sigma$ such that $d_{BM}(K_x,K_y)<1+\delta$ for all $x,y$, and assume that $d_{BM}(K_x, B^2)\geq 1+\epsilon$, for all $x$. We may assume that $\cent(K_x)=0$ for all $x$.
	
	\textit{Step 1: choose parameters $\alpha_3, \beta_0$ such that $G_{\alpha_3}(K_x, K_y)$ does not contain a pair of elements $g_1, g_2$ with $d(g_1, g_2)\in [\beta_0/2, \beta_0]$.}
	
	Denote $\alpha_0=\alpha(\epsilon)$, $\beta_0'=\beta(\epsilon)$ as in eq. \eqref{eq:alpha_beta}. 
	Let $\beta_0\leq \beta_0'$ be the largest real of the form $\beta_0=\frac{2\pi}{B_0}$ with $B_0\geq 3$ an integer. 
	By Proposition \ref{prop:dense_euclidean}, $G^+_{\alpha_0}(K_x)$ is not a $\beta_0$-net. Note that $ \alpha_0,\beta_0\in[c\epsilon, C\epsilon]$.
	
	Denote $\alpha_1=\alpha_0/(2B_0)$, so that $c\epsilon^2\leq \alpha_1\leq C\epsilon^2$. It follows that $G^+_{\alpha_1}(K_x)$ does not contain elements $g$ with  $d(g, \id)\in [\beta_0/2, \beta_0]$: if $g$ were such an element, then the iterates $g^j\in G^+_{\alpha_0}(K_x)$ for $j=0, 1, \dots,  2B_0$ form a $\beta_0$-net, in contradiction. 
	
	Denote $\alpha_2=\alpha_1/2$, $\beta_1=\beta_0/2$. It follows that the distance between $g_1, g_2\in G_{\alpha_2}(K_x)$ cannot lie in $[\beta_1, \beta_0]$, or else
	$g_1^{-1}g_2\in G^+_{\alpha_1}(K_x)$, while $d(g_1^{-1}g_2, \id)\in [\beta_1, \beta_0]$, in contradiction.

	Using Lemma \ref{lem:equivalent_norms}, we let $\delta'=b_2'\delta$ be the upper bound on $d_{BL}(K, L)$ when $d_{BM}(K, L)\leq 1+ \delta$. Assume \begin{equation}\label{eq:mysterious_bound}\alpha_3:=\alpha_2-\delta'>8B_0\delta'\iff \delta'<\frac{1}{8B_0+1}\alpha_2,\end{equation} which holds when $\delta<C\epsilon^3$, as we will henceforth assume to arrive at a contradiction. The reason for the choices in eq. \eqref{eq:mysterious_bound} becomes evident at step 2 of the proof; presently we observe that $\alpha_3>\delta'$, and so $G_{\alpha_3}(K_x, K_y)$  is non-empty.  Furthermore, $G_{\alpha_3}(K_x, K_y)$ does not contain a pair of elements $g_1, g_2$ with 
	$d(g_1, g_2)\in [\beta_1, \beta_0]$: fixing $h\in \OO(\calE(K_x), \calE(K_y))$ such that $d_{BL}(K_x, K_y)=d_{\calE(K_y)}(h K_x, K_y)$, we would otherwise have elements
	$\hat g_1:=g_1h^{-1}, \hat g_2:=g_2h^{-1}\in G_{\alpha_3+\delta'}(K_y)= G_{\alpha_2}(K_y)$, with $d(\hat g_1, \hat g_2)\in  [\beta_1, \beta_0]$, in contradiction.
	
	\textit{Step 2: find a range of parameters $\alpha$ for a fixed point $x_0\in \Sigma$ for which $G_{\alpha}(K_{x_0})$ has the same $\beta_1$-clusters.} 
	
	By a \emph{$\beta$-cluster} in $A\subset S^1$  we understand a maximal subset of $A$ of diameter at most $\beta$. 
	
	Recall that $G_{\alpha_3}(K_x, K_y)\subset\OO(\calE(K_x), \calE(K_y))$, and the latter is the disjoint union of two circles equipped with a metric.

	Let $0<\alpha\leq \alpha_3$. We say that $g_1,g_2\in G_{\alpha}(K_x, K_y)$ are $\beta_1$-close if $d(g_1,g_2)<\beta_1$. By construction, this is an equivalence relation on  $G_{\alpha}(K_x, K_y)$. The $\beta_1$-close equivalence classes are precisely the $\beta_1$-clusters of $G_{\alpha}(K_x, K_y)$, which are therefore pairwise disjoint. Moreover, any two $\beta_1$-clusters  in $G_{\alpha}(K_x, K_y)$ are separated by a distance of at least $\beta_0$. 
	
	Given $g\in G_{\alpha}(K_x, K_y)$, we write $[g]$ for the unique $\beta_1$-cluster containing $g$. Furthermore, for $\alpha\leq \alpha_3$ and $0<\eta<\alpha$, if $Z\subset G_{\alpha-\eta}(K_x, K_y)$ is a $\beta_1$-cluster, and $h\in G_{\eta}(K_y, K_z)$, write $h[Z]:=[h\circ g]\subset  G_{\alpha}(K_x, K_z)$, where $g\in Z$ can be arbitrary.
	
	The number of $\beta_1$-clusters in $G_\alpha(K_x, K_y)$ is a non-decreasing function of $\alpha\in (0, \alpha_3]$, which is bounded from above by $2B_0$.
	Given $0<\alpha<\alpha'\leq \alpha_3$, there is a natural injective map from the $\beta_1$-clusters of $G_\alpha(K_x, K_y)$ into the $\beta_1$-clusters of  $G_{\alpha'}(K_x, K_y)$, given by set inclusion.
	
	Fix $x_0\in \Sigma$. Since $\alpha_3>8B_0\delta'$, we can find $\gamma_0\in [0, \alpha_3-3\delta']$ such that the number of $\beta_1$-clusters in $G_\alpha(K_{x_0})$ is constant for $\alpha\in [\gamma_0, \gamma_0+3\delta']$. Putting $\gamma_1=\gamma_0+3\delta'$, the $\beta_1$-clusters of $G_\alpha(K_{x_0})$ are therefore naturally identified across the range $\alpha\in[\gamma_0, \gamma_1]$.
	
	\textit{Step 3: construct a covering space of $\Sigma$ out of the $\beta_1$-clusters in $G_{\gamma_0}(K_{x_0})$.}
	
	For all $x\in \Sigma$, choose $h_x\in \OO(\calE(K_{x_0}), \calE(K_x))$ such that $d_{\calE(K_x)}(h_x K_{x_0}, K_x)<\frac 32\delta'$. The map $Z\mapsto h_x[Z]$ gives an embedding of the set of $\beta_1$-clusters in $G_{\gamma_0}(K_{x_0})$ into the set of $\beta_1$-clusters in $G_{\gamma_1}(K_{x_0}, K_x)$. Furthermore, the image $C_x$ of this embedding is independent of the choice of $h_x$: if $\hat h_x$ is another such map, and $g\in G_{\gamma_0}(K_{x_0})$, then $g'=\hat h_x^{-1}h_x g\in G_{\gamma_1}(K_{x_0})$. Recall by construction that we may find $g''\in G_{\gamma_0}(K_{x_0})$ with $d(g', g'')<\beta_1$. It then follows that $d(h_x g, \hat h_x g'')=d(\hat h_x g', \hat h_x g'')<\beta_1$, and so $\hat h_x g''$ and $h_x g$ define the same $\beta_1$-cluster in $G_{\gamma_1}(K_{x_0}, K_x)$, proving the claim. 
	Let $m$ be the number of elements (clusters) in $C_x$, which is the same for all $x$ by construction.
	
	Denote by $p:\widehat C\to \Sigma$ the bundle with fiber $C_x$ over $x\in \Sigma$. Our next goal will be topologize $\widehat C$ so as to make it a covering space of $\Sigma$. We will next describe a natural identification of the fibers of $p$ over small subsets of $\Sigma$. 
	
	Fix any Riemannian metric on $\Sigma$, and denote by $B_r(z)$ an open ball of radius $r$ centered at $z$ in this metric. Take $r_0$ small enough so that $B_{r_0}(z)$ can be identified with an open subset of $T_z\Sigma$ through the exponential map, for all $z\in \Sigma$. For $x,y\in B_{r_0}(z)$, $x=\exp_z(\xi)$, $y=\exp_z(\eta)$, we let $\overline u_{xy}^z:T_z\Sigma\to T_z\Sigma$ be the unique positive definite map, with respect to the inner product determined by $d_\xi\exp^{-1}(\calE(K_x))$, mapping $d_\xi\exp ^{-1}(\calE(K_x))$ to $d_\eta \exp^{-1}(\calE(K_y))$. 
	
	Thus $u_{xy}^z:=  d_\eta\exp\circ \overline u_{xy}^z\circ d_\xi\exp ^{-1}\in \SO(\calE(K_x), \calE(K_y))$. Decreasing $r_0$ further, we can achieve $u_{xy}^z\in G^+_{\delta'/4}(K_x, K_y)$ for all $x,y\in B_{r_0}(z)$ by the continuity of $x\mapsto K_x$, as well as $d(u_{xy}^z, u_{xy}^{z'})<\frac\epsilon 2$ for $x,y\in B_{r_0}(z)\cap B_{r_0}(z')$.
	
	Choose $h_x\in G_{\delta'}(K_{x_0}, K_x)$. Now if $Z$ is a $\beta_1$-cluster in $G_{\gamma_0}(K_{x_0})$, then $h_x[Z]$ is a $\beta_1$-cluster in  $C_x$, while  $u_{xy}^z\circ h_x[Z]$ is a $\beta_1$-cluster in $C_y$, and we define the bijection $b_{xy}:p^{-1}(x)\to p^{-1}(y)$ by setting $b_{xy}(h_x[Z])=u^z_{xy}\circ h_x[Z]$, for all $Z$. One readily checks that $b_{xy}$ is independent of the choice of $h_x$, that $b_{xy}$ does not depend on $z$ as long as $x,y\in B_{r_0}(z)$, and that $b_{yx'}\circ b_{xy}=b_{xx'}$ for all $x, y, x'\in B_{r_0}(z)$. 
	
	This now allows us to topologize $\widehat C$: Taking $U\subset B_{r_0}(z)$ open, $p^{-1}(U)$ is, by the above, the disjoint union of $m$ copies $U_1,\dots, U_m$ of $U$. The collection of all such $U_j$ is a basis of the topology of $\widehat C$.  With this topology, $p:\widehat C\to \Sigma$ is an $m$-sheeted covering map. Moreover, $\widehat C$ naturally inherits a smooth structure from $\Sigma$ such that $p$ is a local diffeomorphism.
		
	The bundle $p^*\widehat C$ over $\widehat C$ trivially admits a section $y\mapsto J_y$. Recall that $J_y\subset G_{\gamma_1}(K_{x_0}, K_{p(y)})$. 
	 
	\textit{Step 4: find a continuous section.} 
	We next choose a continuous field of closed intervals $\widehat C \ni y \mapsto I_y\subset \OO(\calE(K_{x_0}), \calE(K_{p(y)}))$ of length $3\beta_1$ such that $J_y\subset I_y$. This can be done as follows. Choose a finite cover of $\widehat C$ by $r_0$-balls $U_j=B_{r_0}(z_j)$. Let $\rho_j$ be a subordinate partition of unity.
	In each $U_j$, one can easily choose a field of  $3\beta_1$-intervals $ I^j_y$ for $y\in U_j$ such that $J_y\subset  I^j_y$, and $u^{z_j}_{p(y)p(z)}( I^j_y)\subset  \OO(\calE(K_{x_0}), \calE(K_{p(z)}))$ is independent of $y$, if $r_0$ is sufficiently small. 
	Then define $I_y:=\sum_j \rho_j(y)  I^j_y$, namely the weighted mean of the intervals $ I^j_y$ with weights $\rho_j(y)$, which is well defined since all of the intervals $I^j_y$ are contained in a fixed interval $I_y'$ of length $2\beta_1+\beta_1+2\beta_1=5\beta_1<2\pi$ containing $J_y$, and so the averaging is carried out inside $I'_y$. Define $P_y$ as the center of $I_y$.

	Thus $P_y$ is a continuous global section over $\widehat C$ of the principal $\OO(2)$-bundle with fiber $\OO(\calE(K_{x_0}), \calE(K_{p(y)}))\subset \GL(T_{x_0}\Sigma, T_y\widehat C)$ over $y\in \widehat C$, which readily implies that $\widehat C$ admits a continuous field of tangent lines: taking any tangent line $L\subset T_{x_0}\Sigma$, $P_yL$ defines such a field. However $\chi(\widehat C)=m\chi(\Sigma)\neq 0$, a contradiction.
	
	This contradiction means that one cannot have $\delta<C\epsilon^3$, concluding the proof.
	
\endproof

\textit{Proof of Corollary \ref{cor:monochromatic}.}
Denote by $K_x\subset T_x\Sigma$ the unit ball of the Finsler structure. By Theorem \ref{thm:main_g}, $d_{BM}(K_x, B^2)<1+c\delta^{1/3}$ for all $x\in \Sigma$.

The assignment $x\mapsto \cent(K_x)$ defines a vector field over $\Sigma$, which must vanish somewhere by the Poincar\'e-Hopf theorem since by assumption $\chi(\Sigma)\neq 0$. Taking $z$ to be a point where $\cent(K_z)=0$, we deduce by Corollary \ref{cor:linear_BM} that $d_{BM}^{\Lin}(K_z, B^2)<1+c\delta^{1/3}$. As $d_{BM}^{\Lin}(K_x, K_z)<1+\delta$ for all $x$, the conclusion follows.
\qed

\section{The stable Banach conjecture for planar sections}\label{sec:banach}

\textit{Proof of Theorem \ref{thm:main_norm}.} 
Let $K$ be the unit ball of $V$. Let $F\subset V$ be a $3$-dimensional subspace.
By Theorem \ref{thm:main_g}, $d_{BM}(K\cap E, B^2)\leq 1+c\delta^{1/3}$ for all planes $E\subset F$

Fix a Euclidean structure on $F$, and denote by $S(F)$ the unit sphere. For a $2$-dimensional subspace $\theta^\perp\subset F$ with $\theta\in S(F)$, define $p(\theta):=\cent(K\cap \theta^\perp)\in \theta^\perp$.
As $p(\theta)$ defines a continuous vector field over $S(F)$, it must vanish somewhere. Let $E_0\subset F$ be a plane such that $\cent(K\cap E_0)=0$.

By Corollary \ref{cor:linear_BM}, there is an absolute constant $c>0$ such that 
\[(1+c\delta^{1/3})^{-1}\calE(K\cap E_0)\subset K\cap E_0\subset (1+c\delta^{1/3})\calE(K\cap E_0).\]
In particular, $K\cap E_0$ is nearly centrally-symmetric: \[-(K\cap E_0)\subset (1+c\delta^{1/3})(K\cap E_0).\]
By assumption, it holds for any other linear plane $E\subset V$ that $d^{\Lin}_{BM}(K\cap E, K\cap E_0)<1+\delta$, and it easily follows that
$-(K\cap E)\subset (1+c\delta^{1/3})(K\cap E)$ for all $E$. Therefore, $-K\subset (1+c\delta^{1/3})K$.

 Let $K_s=\frac{1}{2}(K+(-K))$ be the Minkowski symmetrization of $K$. It follows readily that
 \begin{equation}\label{eq:symmetric_equiv} (1+c\delta^{1/3})^{-1}K_s\subset K\subset (1+c\delta^{1/3})K_s.\end{equation}

 It follows by the above that $d_{BM}(K_s\cap E, B^2)<1+c\delta^{1/3}$ for all linear planes $E$. By Lemma \ref{lem:katoetal}, the von Neumann--Jordan constant of each section $K_s\cap E$ satisfies $C_{NJ}(K_s\cap E)< 1+c\delta^{1/3}$. Consequently, $C_{NJ}(K_s) <1+ c\delta^{1/3}$. 
By Theorem \ref{thm:passer}, 
\begin{align*}d_{BM}(K_s, B^n)=d^{\Lin}_{BM}(K_s, B^n) <1+cn^2 \delta^{1/3},\end{align*}
and by eq. \eqref{eq:symmetric_equiv}, $d_{BM}^{\Lin}(K, B^n)<1+cn^2 \delta^{1/3}$.
\qed\endproof

To deduce Theorem \ref{thm:main1}, we will replace Passer's theorem with Theorem \ref{mprop:one_center}, which we now proceed to prove and restate here for convenience. It is a stable version of the simple geometric fact proved in \cite{AuerbachMazurUlam35} (see also \cite[Lemma 16.12]{busemann_geodesics}), asserting that if the $2$-dimensional sections of a convex body through a fixed interior point are all ellipses, then the body is an ellipsoid. 

\begin{Theorem}\label{prop:one_center}
	Let $K\subset \R^n$ be a convex body with $0\in\mathrm{int}(K)$ such that for all $2$-dimensional linear planes $E$, $d_{BM}(K\cap E, B^2)<1+\epsilon$.
	Then $d_{BM}(K, B^n)<1+C_n\sqrt \epsilon$, where $C_n=cn^{2n^2}$.
\end{Theorem}
\proof
We prove the statement with $C'_n=c^nd_{n-1}d_{n-2}\cdots d_2$ replacing $C_n$, as $C'_n<C_n$. We proceed by induction on $n$, with the base $n=2$ being trivial.
	
	First consider the map $\Gr_{n-1}(\R^n)\ni H\mapsto \cent(K\cap H)\in H$. This is a global section of the vector bundle $\gamma^\perp$ with fiber $H$ over $\mathbb R\mathbb P^{n-1}=\mathbb P((\R^n)^*)$, which has full Stiefel-Whitney class $w(\gamma^\perp)=1+a+\dots+a^{n-1}$, where $a\in H^1(\mathbb R\mathbb P^{n-1}, \mathbb Z_2)$ is the generator (see e.g. \cite[chapter 4]{milnor_stasheff}). It follows that any global section must vanish somewhere, and we fix a hyperplane $H_0$ such that $\cent(K\cap H_0)=0$. Denote $\mathcal E_0=\calE(K\cap H_0)$. 
	
	By the induction hypothesis, $d_{BM}(K\cap H_0, B^{n-1})<1+C'_{n-1}\sqrt{\epsilon}$. By Corollary \ref{cor:linear_BM}, 
	\begin{equation}\label{eq:hyperplane_ellipsoid}(1+d_{n-1}C'_{n-1}\sqrt\epsilon)^{-1}\calE_0\subset K\cap H_0\subset(1+d_{n-1}C'_{n-1}\sqrt\epsilon)\calE_0.\end{equation} Since any chord of $\calE_0$ through the origin is bisected by it, it follows that if $[A, A']$ is a chord of $K\cap H_0$ through the origin, then for any Euclidean norm one has
	\begin{equation}\label{eq:almost_middle}(1+d_{n-1}C'_{n-1}\sqrt\epsilon)^{-1}\leq\frac{|A|}{|A'|}\leq 1+d_{n-1}C'_{n-1}\sqrt\epsilon\end{equation}
	
	Let $H_1, H_2$ be the affine hyperplanes supporting $K$ and parallel to $H_0$. Assume $H_1$ is at least as far from the origin as $H_2$ (with respect to any norm on $\R^n/H_0$). Fix $z\in \partial K\cap H_1$. 
	\begin{center}

\tikzset {_ks93ug22d/.code = {\pgfsetadditionalshadetransform{ \pgftransformshift{\pgfpoint{0 bp } { 0 bp }  }  \pgftransformrotate{0 }  \pgftransformscale{2 }  }}}
\pgfdeclarehorizontalshading{_wxdrotob3}{150bp}{rgb(0bp)=(1,1,1);
rgb(37.5bp)=(1,1,1);
rgb(62.5bp)=(0.9,0.9,0.9);
rgb(100bp)=(0.9,0.9,0.9)}
\tikzset{every picture/.style={line width=0.75pt}} 
 \begin{figure} 
\begin{tikzpicture}[x=0.75pt,y=0.75pt,yscale=-1,xscale=1]

\path  [shading=_wxdrotob3,_ks93ug22d] (174.58,1259.52) -- (246.45,1280.51) -- (347.68,1344.17) -- (173.04,1412.57) -- (98.08,1361.78) -- (87.26,1300.83) -- cycle ; 
 \draw   (174.58,1259.52) -- (246.45,1280.51) -- (347.68,1344.17) -- (173.04,1412.57) -- (98.08,1361.78) -- (87.26,1300.83) -- cycle ; 

\draw  [fill={rgb, 255:red, 0; green, 0; blue, 0 }  ,fill opacity=1 ] (179.6,1380.73) .. controls (179.6,1379.8) and (180.47,1379.04) .. (181.54,1379.04) .. controls (182.6,1379.04) and (183.47,1379.8) .. (183.47,1380.73) .. controls (183.47,1381.67) and (182.6,1382.43) .. (181.54,1382.43) .. controls (180.47,1382.43) and (179.6,1381.67) .. (179.6,1380.73) -- cycle ;
\draw    (109.28,1369.91) -- (9.6,1394.96) ;
\draw    (9.6,1394.96) -- (340.33,1395.64) ;
\draw  [dash pattern={on 0.84pt off 2.51pt}]  (109.28,1369.91) -- (281.61,1369.91) ;
\draw    (281.61,1369.91) -- (397.52,1369.23) ;
\draw    (340.33,1395.64) -- (397.52,1369.23) ;
\draw  [fill={rgb, 255:red, 0; green, 0; blue, 0 }  ,fill opacity=1 ] (170.33,1412.23) .. controls (170.33,1411.3) and (171.2,1410.54) .. (172.26,1410.54) .. controls (173.33,1410.54) and (174.19,1411.3) .. (174.19,1412.23) .. controls (174.19,1413.17) and (173.33,1413.93) .. (172.26,1413.93) .. controls (171.2,1413.93) and (170.33,1413.17) .. (170.33,1412.23) -- cycle ;
\draw  [fill={rgb, 255:red, 0; green, 0; blue, 0 }  ,fill opacity=1 ] (172.65,1259.52) .. controls (172.65,1258.58) and (173.51,1257.83) .. (174.58,1257.83) .. controls (175.65,1257.83) and (176.51,1258.58) .. (176.51,1259.52) .. controls (176.51,1260.45) and (175.65,1261.21) .. (174.58,1261.21) .. controls (173.51,1261.21) and (172.65,1260.45) .. (172.65,1259.52) -- cycle ;
\draw   (120.1,1247.67) -- (392.88,1247.67) -- (311.74,1272.72) -- (38.96,1272.72) -- cycle ;
\draw    (11.92,1427.13) -- (310.97,1426.12) ;
\draw    (310.97,1426.12) -- (368.15,1399.7) ;
\draw    (11.92,1427.13) -- (94.6,1400.72) ;
\draw    (154.88,1400.72) -- (94.6,1400.72) ;
\draw  [dash pattern={on 0.84pt off 2.51pt}]  (154.88,1400.72) -- (205.1,1400.04) ;
\draw    (205.1,1400.04) -- (368.15,1399.7) ;
\draw    (172.6,1216.8) -- (183.34,1405.69) -- (186.6,1475.8) ;
\draw  [dash pattern={on 4.5pt off 4.5pt}]  (83.8,1387.22) -- (180.97,1380.67) -- (264.62,1375.03) ;
\draw  [dash pattern={on 4.5pt off 4.5pt}]  (110.84,1265.99) -- (290.89,1252.45) ;
\draw    (107.75,1272.77) -- (83.8,1387.22) ;
\draw    (62.16,1481.35) -- (241.42,1473.23) ;
\draw    (81.47,1394.66) -- (75.3,1421.76) ;
\draw    (120.9,1219.24) -- (110.84,1265.99) ;
\draw    (120.9,1219.24) -- (300.15,1211.12) ;
\draw    (300.15,1211.12) -- (291.76,1250.42) ;
\draw    (289.35,1272.09) -- (281.62,1303.92) ;
\draw  [dash pattern={on 0.84pt off 2.51pt}]  (281.62,1303.92) -- (264.62,1375.03) ;
\draw    (263.06,1378.4) -- (254.57,1410.92) ;
\draw    (73.74,1427.16) -- (62.16,1481.35) ;
\draw    (241.42,1473.23) -- (252.26,1425.82) ;
\draw  [dash pattern={on 4.5pt off 4.5pt}]  (75.3,1421.76) -- (254.57,1410.92) ;

\draw (399.09,1359.23) node [anchor=north west][inner sep=0.75pt]    {$H_{0}$};
\draw (398.31,1238.68) node [anchor=north west][inner sep=0.75pt]    {$H_{1}$};
\draw (376.38,1397.03) node [anchor=north west][inner sep=0.75pt]    {$H_{2}$};
\draw (159.02,1251.37) node [anchor=north west][inner sep=0.75pt]    {$z$};
\draw (186.87,1381.6) node [anchor=north west][inner sep=0.75pt]    {$o$};
\draw (270.3,1218.38) node [anchor=north west][inner sep=0.75pt]    {$E$};
\draw (92.38,1372.24) node [anchor=north west][inner sep=0.75pt]    {$L$};
\draw (120.1,1249.67) node [anchor=north west][inner sep=0.75pt]    {$L_{1}$};
\draw (85.62,1404.37) node [anchor=north west][inner sep=0.75pt]    {$L_{2}$};
\draw (303.34,1328.65) node [anchor=north west][inner sep=0.75pt]    {$K$};

\end{tikzpicture}
\caption{}
\end{figure}
	\end{center}
	Now let $E$ be any linear $2$-dimensional plane through $z$, and define $L=E\cap H_0$. Denote $\calE_E:=\calE(K\cap E)$, $p_E:=\cent(\calE_E)=\cent(K\cap E)$. We will use the Euclidean structure on $E$ with unit ball $\calE_E-p_E$. Note that by Corollary \ref{cor:linear_BM}, \begin{equation}\label{eq:E_ellipse}(1+c\epsilon)^{-1}(\calE_E-p_E)\subset K\cap E-p_E\subset(1+c\epsilon)(\calE_E-p_E).\end{equation}
	
	Let $L_i=H_i\cap E$, $i=1,2$ be the lines in $E$ which are supporting lines of $K\cap E$ and are parallel to $L$. As the distance between $L_1$ and $L_2$ is at least the diameter of $(1+c\epsilon)^{-1}(\calE_E-p_E)$, it follows that \begin{equation}\label{eq:larger_distance} \mathrm{dist}(L_1, 0)\geq 1-c\epsilon.\end{equation}
	
	We will repeatedly make use of the following simple fact: If two concentric discs $B_1, B_2$ of radii $1-c\epsilon$, resp. $1+C\epsilon$ are given, and $L$ is a line, then for $\epsilon>0$, we have \begin{equation}\label{eq:concentric_line}\mathrm{Length}(L\cap (B_2\setminus B_1))< c\sqrt \epsilon.\end{equation}

	Orient $L$ arbitrarily. Denote $[A, A']=K\cap L$, $[A_E, A'_E]=\calE_E\cap L$.  By eq. \eqref{eq:almost_middle}, 
	\[ |\frac{A+A'}{2}-0|\leq cd_{n-1}C'_{n-1}\sqrt\epsilon.\]
	By eqs. \eqref{eq:E_ellipse} and \eqref{eq:concentric_line}, $|A-A_E|\leq c\sqrt{\epsilon}$, and $|A'-A_E'|\leq c\sqrt{\epsilon}$. Consequently, setting $0_E:=(A_E+A'_E)/2$, we find
	\begin{equation}\label{eq:origins_close}|0_E-0|\leq cd_{n-1}C'_{n-1}\sqrt{\epsilon}.\end{equation}
	
	Let $z_E$ be the touching point on $\partial \calE_E$ of a translate of $L$ that is on the same side of $L$ as $L_1$. By eq. \eqref{eq:concentric_line}, 
	\begin{equation}\label{eq:z_close}|z-z_E|\leq c\sqrt\epsilon.\end{equation}
		
	Let $[z,w]$ be the chord in $K\cap E$ defined by the line through $z, 0$. Let $[z_E, w_E]$ be the chord in $\calE_E$ defined by the line through $z_E, 0_E$.
	Let $p$ be the midpoint of $[z,w]$, and note that the midpoint of $[z_E, w_E]$ is $p_E$. 
	\vspace{2mm}

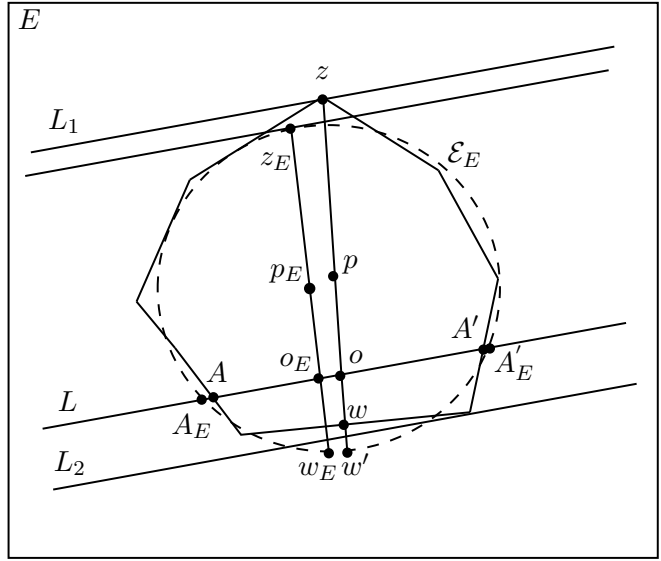
\begin{figure}
\tikzset{every picture/.style={line width=0.75pt}} 

\begin{tikzpicture}[x=0.75pt,y=0.75pt,yscale=-1,xscale=1]

\draw    (473.88,2556.35) -- (532.41,2593.34) ;
\draw    (532.41,2593.34) -- (562.5,2648.05) ;
\draw    (562.5,2648.05) -- (548.29,2715.45) ;
\draw    (399.51,2682.94) -- (432.84,2726.82) -- (548.29,2715.45) ;
\draw    (380.14,2659.6) -- (399.51,2682.94) ;
\draw    (407.03,2598.1) -- (380.14,2659.6) ;
\draw    (407.03,2598.1) -- (473.88,2556.35) ;
\draw   (315.6,2508.8) -- (643.02,2508.8) -- (643.02,2788.96) -- (315.6,2788.96) -- cycle ;
\draw  [fill={rgb, 255:red, 0; green, 0; blue, 0 }  ,fill opacity=1 ] (480.78,2697.06) .. controls (480.78,2695.9) and (481.71,2694.96) .. (482.86,2694.96) .. controls (484.01,2694.96) and (484.94,2695.9) .. (484.94,2697.06) .. controls (484.94,2698.23) and (484.01,2699.17) .. (482.86,2699.17) .. controls (481.71,2699.17) and (480.78,2698.23) .. (480.78,2697.06) -- cycle ;
\draw  [fill={rgb, 255:red, 0; green, 0; blue, 0 }  ,fill opacity=1 ] (471.9,2557.59) .. controls (471.9,2556.43) and (472.83,2555.48) .. (473.98,2555.48) .. controls (475.13,2555.48) and (476.06,2556.43) .. (476.06,2557.59) .. controls (476.06,2558.75) and (475.13,2559.69) .. (473.98,2559.69) .. controls (472.83,2559.69) and (471.9,2558.75) .. (471.9,2557.59) -- cycle ;
\draw    (332.77,2723.76) -- (627.03,2670.37) ;
\draw    (326.85,2584.28) -- (621.11,2530.9) ;
\draw    (338.1,2754.47) -- (632.36,2701.08) ;
\draw  [dash pattern={on 4.5pt off 4.5pt}] (390.8,2652.88) .. controls (390.8,2607.39) and (429.45,2570.51) .. (477.12,2570.51) .. controls (524.79,2570.51) and (563.43,2607.39) .. (563.43,2652.88) .. controls (563.43,2698.36) and (524.79,2735.24) .. (477.12,2735.24) .. controls (429.45,2735.24) and (390.8,2698.36) .. (390.8,2652.88) -- cycle ;
\draw  [fill={rgb, 255:red, 0; green, 0; blue, 0 }  ,fill opacity=1 ] (465.1,2652.99) .. controls (465.1,2651.66) and (466.16,2650.58) .. (467.47,2650.58) .. controls (468.78,2650.58) and (469.84,2651.66) .. (469.84,2652.99) .. controls (469.84,2654.31) and (468.78,2655.39) .. (467.47,2655.39) .. controls (466.16,2655.39) and (465.1,2654.31) .. (465.1,2652.99) -- cycle ;
\draw  [fill={rgb, 255:red, 0; green, 0; blue, 0 }  ,fill opacity=1 ] (416.84,2707.88) .. controls (416.84,2706.72) and (417.77,2705.78) .. (418.92,2705.78) .. controls (420.07,2705.78) and (421,2706.72) .. (421,2707.88) .. controls (421,2709.05) and (420.07,2709.99) .. (418.92,2709.99) .. controls (417.77,2709.99) and (416.84,2709.05) .. (416.84,2707.88) -- cycle ;
\draw  [fill={rgb, 255:red, 0; green, 0; blue, 0 }  ,fill opacity=1 ] (410.92,2709.09) .. controls (410.92,2707.92) and (411.85,2706.98) .. (413,2706.98) .. controls (414.15,2706.98) and (415.08,2707.92) .. (415.08,2709.09) .. controls (415.08,2710.25) and (414.15,2711.19) .. (413,2711.19) .. controls (411.85,2711.19) and (410.92,2710.25) .. (410.92,2709.09) -- cycle ;
\draw  [fill={rgb, 255:red, 0; green, 0; blue, 0 }  ,fill opacity=1 ] (553.01,2683.84) .. controls (553.01,2682.68) and (553.94,2681.73) .. (555.09,2681.73) .. controls (556.24,2681.73) and (557.17,2682.68) .. (557.17,2683.84) .. controls (557.17,2685) and (556.24,2685.94) .. (555.09,2685.94) .. controls (553.94,2685.94) and (553.01,2685) .. (553.01,2683.84) -- cycle ;
\draw  [fill={rgb, 255:red, 0; green, 0; blue, 0 }  ,fill opacity=1 ] (484.39,2735.84) .. controls (484.39,2734.68) and (485.32,2733.74) .. (486.47,2733.74) .. controls (487.61,2733.74) and (488.54,2734.68) .. (488.54,2735.84) .. controls (488.54,2737) and (487.61,2737.94) .. (486.47,2737.94) .. controls (485.32,2737.94) and (484.39,2737) .. (484.39,2735.84) -- cycle ;
\draw  [fill={rgb, 255:red, 0; green, 0; blue, 0 }  ,fill opacity=1 ] (556.28,2683.26) .. controls (556.28,2682.09) and (557.21,2681.15) .. (558.35,2681.15) .. controls (559.5,2681.15) and (560.43,2682.09) .. (560.43,2683.26) .. controls (560.43,2684.42) and (559.5,2685.36) .. (558.35,2685.36) .. controls (557.21,2685.36) and (556.28,2684.42) .. (556.28,2683.26) -- cycle ;
\draw  [fill={rgb, 255:red, 0; green, 0; blue, 0 }  ,fill opacity=1 ] (469.84,2698.37) .. controls (469.84,2697.21) and (470.77,2696.27) .. (471.91,2696.27) .. controls (473.06,2696.27) and (473.99,2697.21) .. (473.99,2698.37) .. controls (473.99,2699.54) and (473.06,2700.48) .. (471.91,2700.48) .. controls (470.77,2700.48) and (469.84,2699.54) .. (469.84,2698.37) -- cycle ;
\draw    (474.28,2556.92) -- (483.39,2695.18) -- (486.47,2735.84) ;
\draw  [fill={rgb, 255:red, 0; green, 0; blue, 0 }  ,fill opacity=1 ] (482.56,2721.78) .. controls (482.56,2720.62) and (483.49,2719.68) .. (484.64,2719.68) .. controls (485.78,2719.68) and (486.71,2720.62) .. (486.71,2721.78) .. controls (486.71,2722.95) and (485.78,2723.89) .. (484.64,2723.89) .. controls (483.49,2723.89) and (482.56,2722.95) .. (482.56,2721.78) -- cycle ;
\draw  [fill={rgb, 255:red, 0; green, 0; blue, 0 }  ,fill opacity=1 ] (455.67,2572.32) .. controls (455.67,2571.16) and (456.6,2570.21) .. (457.75,2570.21) .. controls (458.89,2570.21) and (459.82,2571.16) .. (459.82,2572.32) .. controls (459.82,2573.48) and (458.89,2574.42) .. (457.75,2574.42) .. controls (456.6,2574.42) and (455.67,2573.48) .. (455.67,2572.32) -- cycle ;
\draw    (458.05,2571.72) -- (477.12,2735.24) ;
\draw  [fill={rgb, 255:red, 0; green, 0; blue, 0 }  ,fill opacity=1 ] (475.04,2736.14) .. controls (475.04,2734.98) and (475.97,2734.04) .. (477.12,2734.04) .. controls (478.27,2734.04) and (479.2,2734.98) .. (479.2,2736.14) .. controls (479.2,2737.3) and (478.27,2738.25) .. (477.12,2738.25) .. controls (475.97,2738.25) and (475.04,2737.3) .. (475.04,2736.14) -- cycle ;
\draw  [fill={rgb, 255:red, 0; green, 0; blue, 0 }  ,fill opacity=1 ] (477.23,2646.78) .. controls (477.23,2645.62) and (478.16,2644.67) .. (479.31,2644.67) .. controls (480.46,2644.67) and (481.39,2645.62) .. (481.39,2646.78) .. controls (481.39,2647.94) and (480.46,2648.88) .. (479.31,2648.88) .. controls (478.16,2648.88) and (477.23,2647.94) .. (477.23,2646.78) -- cycle ;

\draw (338.98,2703.3) node [anchor=north west][inner sep=0.75pt]    {$L$};
\draw (334,2561.02) node [anchor=north west][inner sep=0.75pt]    {$L_{1}$};
\draw (336.96,2733.56) node [anchor=north west][inner sep=0.75pt]    {$L_{2}$};
\draw (468.6,2539.37) node [anchor=north west][inner sep=0.75pt]    {$z$};
\draw (485.66,2682.67) node [anchor=north west][inner sep=0.75pt]    {$o$};
\draw (446.08,2639.17) node [anchor=north west][inner sep=0.75pt]    {$p_{E}$};
\draw (535.48,2577.85) node [anchor=north west][inner sep=0.75pt]    {$\mathcal{E}_{E}$};
\draw (413.56,2688.46) node [anchor=north west][inner sep=0.75pt]    {$A$};
\draw (537.07,2665.02) node [anchor=north west][inner sep=0.75pt]    {$A'$};
\draw (396.33,2714.7) node [anchor=north west][inner sep=0.75pt]    {$A_{E}$};
\draw (558.78,2680.01) node [anchor=north west][inner sep=0.75pt]    {$A_{E}^{'}$};
\draw (450.76,2684.48) node [anchor=north west][inner sep=0.75pt]    {$o_{E}$};
\draw (484.66,2710.57) node [anchor=north west][inner sep=0.75pt]    {$w$};
\draw (440.77,2583.11) node [anchor=north west][inner sep=0.75pt]    {$z_{E}$};
\draw (460.17,2739.53) node [anchor=north west][inner sep=0.75pt]    {$w_{E}$};
\draw (482.89,2634.63) node [anchor=north west][inner sep=0.75pt]    {$p$};
\draw (481.78,2734.56) node [anchor=north west][inner sep=0.75pt]    {$w'$};
\draw (318.85,2510.7) node [anchor=north west][inner sep=0.75pt]    {$E$};
\draw    (323.99,2596.21) -- (618.25,2542.82); 

\end{tikzpicture}	
\caption{Constructions in the plane $E$}
\end{figure}
	\vspace{2mm}
	We next claim that \begin{equation}\label{eq:w_close}|w-w_E|\leq cd_{n-1}C'_{n-1}\sqrt{\epsilon}.\end{equation} Denoting by $w'$ the intersection of the ray $R$ from $z$ to $w$ with $\partial \calE_E$, it holds by eq. \eqref{eq:concentric_line} that $|w-w'|\leq c\sqrt \epsilon$. Next, from eq. \eqref{eq:larger_distance} we conclude that the angle formed by the ray $R_E$ from $z_E$ to $0_E$ with $\partial \calE_E$ lies in $[\frac \pi 4, \frac \pi 2]$. From eqs. \eqref{eq:origins_close}, \eqref{eq:z_close}, and \eqref{eq:larger_distance}
	we deduce that the angle between $R$ and $R_E$ is at most $cd_{n-1}C'_{n-1}\sqrt\epsilon$. As $\mathrm{dist}(R\cap \calE_E, R_E\cap \calE_E)\leq c\sqrt \epsilon$, 
	we conclude that $|w_E-w'|\leq  cd_{n-1}C'_{n-1}\sqrt{\epsilon}$, and so $|w-w_E|\leq |w-w'|+|w'-w_E|\leq   cd_{n-1}C'_{n-1}\sqrt{\epsilon}$, showing eq. \eqref{eq:w_close}.
	
	We deduce from eqs. \eqref{eq:z_close}, \eqref{eq:w_close} that also $|p-p_E|\leq cd_{n-1}C'_{n-1}\sqrt \epsilon$.
	
	From \eqref{eq:E_ellipse} it follows that $K\cap E$ is close  to an ellipse with a center at $p$: 
	\begin{equation}\label{eq:E_ellipse_approx} (1+cd_{n-1}C'_{n-1}\sqrt \epsilon)^{-1}(\calE_E-p) \subset K\cap E-p \subset  (1+cd_{n-1}C_{n-1}\sqrt \epsilon)(\calE_E-p).\end{equation}
	
	We now let $\calE$ be the unique ellipsoid with center at $p$ such that both $z$ and $\partial \calE_0$ lie on $\partial \calE$.
	It follows from eqs. \eqref{eq:E_ellipse_approx}, \eqref{eq:hyperplane_ellipsoid}, \eqref{eq:z_close}, \eqref{eq:larger_distance} that 
	\[  (1+cd_{n-1}C'_{n-1}\sqrt\epsilon)^{-1}(\calE_E-p)\subset(\calE -p)\cap E\subset (1+cd_{n-1}C'_{n-1}\sqrt\epsilon)(\calE_E-p),\]
	and thus also
	\[(1+cd_{n-1}C'_{n-1}\sqrt \epsilon)^{-1}(\calE-p) \subset K-p \subset  (1+cd_{n-1}C'_{n-1}\sqrt \epsilon)(\calE-p).\]
	Taking $C'_n=cd_{n-1}C'_{n-1}$ completes the induction and the proof. 
	
\endproof

\textit{Proof of Theorem \ref{thm:main1}.}
From Theorem \ref{thm:main_g} it follows that $d_{BM}(K\cap E, B^2)<1+c\delta^{1/3}$ for all linear planes $E$. We then deduce from Theorem \ref{prop:one_center} that $d_{BM}(K, B^n)<1+cn^{2n^2}\delta^{1/6}$, concluding the proof.\qed

\bibliographystyle{amsplain} 
\bibliography{shoulderofgiants}

\providecommand{\bysame}{\leavevmode\hbox to3em{\hrulefill}\thinspace}
\providecommand{\MR}{\relax\ifhmode\unskip\space\fi MR }
\providecommand{\MRhref}[2]{%
  \href{http://www.ams.org/mathscinet-getitem?mr=#1}{#2}
}
\providecommand{\href}[2]{#2}
\begin{thebibliography}{10}

\bibitem{artstein_etal}
Shiri Artstein-Avidan, Apostolos Giannopoulos, and Vitali~D. Milman,
  \emph{Asymptotic geometric analysis. {P}art {I}}, Mathematical Surveys and
  Monographs, vol. 202, American Mathematical Society, Providence, RI, 2015.
  \MR{3331351}

\bibitem{AuerbachMazurUlam35}
H~Auerbach, S~Mazur, and S~Ulam, \emph{On a characteristic property of the
  ellipso{\"\i}de}, Monatshefte f{\"u}r Mathematik und Physik \textbf{42}
  (1935), 45--48.

\bibitem{Banach30}
S~Banach, \emph{Th{\'e}orie des op{\'e}rations lin{\'e}aires. monografie
  matematyczne, 1, polskie towarzystwo matematyczne, warszawa, 1932.[5},
  Th{\'e}or{\`e}me sur les ensembles de premi{\`e}re cat{\'e}gorie. Fund. Math
  \textbf{16} (1930), 395--398.

\bibitem{BorLamonedaDesantiagoMontejano21}
Gil Bor, Luis Hern{\'a}ndez~Lamoneda, Valent{\'\i}n Jim{\'e}nez-Desantiago, and
  Luis Montejano, \emph{On the isometric conjecture of {B}anach}, Geometry \&
  Topology \textbf{25} (2021), no.~5, 2621--2642.

\bibitem{montejano_complex}
Javier Bracho and Luis Montejano, \emph{On the complex {Banach} conjecture}, J.
  Convex Anal. \textbf{28} (2021), no.~4, 1211--1222 (English).

\bibitem{BurgerSchneider93}
Thomas Burger and Rolf Schneider, \emph{On convex bodies close to ellipsoids},
  Journal of Geometry \textbf{47} (1993), 16--22.

\bibitem{busemann_geodesics}
Herbert Busemann, \emph{The geometry of geodesics}, reprint of the 1955
  original ed., Mineola, NY: Dover Publications, 2005 (English).

\bibitem{Groemer94}
H~Groemer, \emph{Stability theorems for ellipsoids and spheres}, Journal of the
  London Mathematical Society \textbf{49} (1994), no.~2, 357--370.

\bibitem{Gromov67}
ML~Gromov, \emph{On a geometric hypothesis of {B}anach}, Izv. Akad. Nauk SSSR
  Ser. Mat \textbf{31} (1967), 1105--1114.

\bibitem{Gruber97}
Peter~M Gruber, \emph{Stability of {B}laschke's characterization of ellipsoids
  and {R}adon norms}, Discrete \& Computational Geometry \textbf{17} (1997),
  411--427.

\bibitem{Ivanov18}
Sergei Ivanov, \emph{Monochromatic {F}insler surfaces and a local ellipsoid
  characterization}, Proceedings of the American Mathematical Society
  \textbf{146} (2018), no.~4, 1741--1755.

\bibitem{IvanovMamaevNordskova23}
Sergei Ivanov, Daniil Mamaev, and Anya Nordskova, \emph{{B}anach’s isometric
  subspace problem in dimension four}, Inventiones mathematicae \textbf{233}
  (2023), no.~3, 1393--1425.

\bibitem{ivanov2023local}
Sergei Ivanov, Daniil Mamaev, and Anya Nordskova, \emph{Local {K}akutani's
  ellipsoid characterization}, 2023.

\bibitem{neumann_jordan}
P.~Jordan and J.~von Neumann, \emph{On inner products in linear, metric
  spaces}, Ann. Math. (2) \textbf{36} (1935), 719--723 (English).

\bibitem{KatoMalingradaTakahashi01}
Mikio Kato, Lech Maligranda, and Yasuji Takahashi, \emph{On {J}ames and
  {J}ordan-von {N}eumann constants and the normal structure coefficient of
  {B}anach spaces}, Studia Mathematica \textbf{144} (2001), no.~3, 275--295.

\bibitem{MatveevTroyanov12}
Vladimir~S Matveev and Marc Troyanov, \emph{The {B}inet--{L}egendre metric in
  {F}insler geometry}, Geometry \& Topology \textbf{16} (2012), no.~4,
  2135--2170.

\bibitem{MatveevTroyanov15}
\bysame, \emph{Completeness and incompleteness of the {B}inet--{L}egendre
  metric}, European Journal of Mathematics \textbf{1} (2015), no.~3, 483--502.

\bibitem{milnor_stasheff}
John~W. Milnor and James~D. Stasheff, \emph{Characteristic classes}, Texts
  Read. Math., vol.~32, New Delhi: Hindustan Book Agency, 2005 (English).

\bibitem{Montejano21arxivsurvey}
Luis Montejano, \emph{Convex bodies all whose sections (projections) are
  equal}, European congress of mathematics. Proceedings of the 8th congress,
  8ECM, Portoro\v{z}, Slovenia, June 20--26, 2021, Berlin: European
  Mathematical Society (EMS), 2023, pp.~857--883 (English).

\bibitem{Passer15}
Benjamin Passer, \emph{An approximate version of the {J}ordan von {N}eumann
  theorem for finite-dimensional real normed spaces}, Linear and Multilinear
  Algebra \textbf{63} (2015), no.~1, 68--77.

\bibitem{Schneider14}
Rolf Schneider, \emph{Convex bodies: the {B}runn--{M}inkowski theory}, no. 151,
  Cambridge university press, 2014.

\end{thebibliography}

 \end{document}